\documentclass[a4,10pt]{article}
\usepackage{amsmath,amscd,amssymb}
\usepackage[margin=.8in]{geometry}
\usepackage{xcolor}
\usepackage{mathdots}

\newcommand{\nbiga}{\mathcal{A}}

\newcommand{\nbigf}{\mathcal{F}}

\newcommand{\nbigh}{\mathcal{H}}

\newcommand{\nbigo}{\mathcal{O}}

\newcommand{\nbigs}{\mathcal{S}}
\newcommand{\nbigt}{\mathcal{T}}

\newcommand{\seisuu}{{\mathbb Z}}

\newcommand{\cnum}{{\mathbb C}}
\newcommand{\real}{{\mathbb R}}

\newcommand{\hyperk}{\mathbb{K}}

\newcommand{\vece}{{\boldsymbol e}}

\newcommand{\vecv}{{\boldsymbol v}}

\newcommand{\lrarr}{\longrightarrow}

\newcommand{\pf}{{\bf Proof}\hspace{.1in}}
\newcommand{\qed}{\mbox{\rule{1.2mm}{3mm}}}

\def\Hom{\mathop{\rm Hom}\nolimits}

\def\End{\mathop{\rm End}\nolimits}

\def\Gr{\mathop{\rm Gr}\nolimits}

\def\rank{\mathop{\rm rank}\nolimits}

\def\Gr{\mathop{\rm Gr}\nolimits}

\def\tr{\mathop{\rm tr}\nolimits}
\def\Tr{\mathop{\rm Tr}\nolimits}

\def\id{\mathop{\rm id}\nolimits}

\def\diag{\mathop{\rm diag}\nolimits}

\newcommand{\del}{\partial}
\newcommand{\delbar}{\overline{\del}}

\newcommand{\barz}{\overline{z}}
\newcommand{\zbar}{\barz}

\newcommand{\Sp}{{\mathcal Sp}}

\def\Harm{\mathop{\rm Harm}\nolimits}

\newcommand{\Etilde}{\widetilde{E}}

\newcommand{\htilde}{\widetilde{h}}

\newcommand{\ftilde}{\widetilde{f}}

\newcommand{\Xtilde}{\widetilde{X}}

\newcommand{\Wtilde}{\widetilde{W}}

\newcommand{\Ytilde}{\widetilde{Y}}

\newcommand{\ctilde}{\widetilde{c}}
\newcommand{\taubar}{\overline{\tau}}

\newcommand{\vecq}{\boldsymbol q}

\newtheorem{thm}{Theorem}[section]
\newtheorem{cor}[thm]{Corollary}
\newtheorem{rem}[thm]{Remark}
\newtheorem{lem}[thm]{Lemma}
\newtheorem{prop}[thm]{Proposition}
\newtheorem{df}[thm]{Definition}

\newtheorem{condition}[thm]{Condition}
\newtheorem{question}[thm]{Question}

\begin{document}
\title{Higgs bundles in the Hitchin section over non-compact hyperbolic surfaces}
\author{Qiongling Li\thanks{Chern Institute of Mathematics and LPMC, Nankai University, Tianjin 300071, China, qiongling.li@nankai.edu.cn}
\and Takuro Mochizuki\thanks{Research Institute for Mathematical Sciences, Kyoto University, Kyoto 606-8512, Japan, takuro@kurims.kyoto-u.ac.jp}}
\date{}

\maketitle
\begin{abstract}
    Let $X$ be an arbitrary non-compact hyperbolic Riemann surface, that is, not $\mathbb C$ or $\mathbb C^*$. Given a tuple of holomorphic differentials $\vecq=(q_2,\cdots,q_n)$ on $X$, one can define a Higgs bundle $(\hyperk_{X,n},\theta(\vecq))$ in the Hitchin section. We show there exists a harmonic metric $h$ on $(\hyperk_{X,n},\theta(\vecq))$ satisfying (i) $h$ weakly dominates $h_X$; (ii) $h$ is compatible with the real structure. Here $h_X$ is the Hermitian metric on $\hyperk_{X,n}$ induced by the conformal complete hyperbolic metric $g_X$ on $X.$ Moreover, when $q_i(i=2,\cdots,n)$ are bounded with respect to $g_X$, we show such a harmonic metric on $(\hyperk_{X,n},\theta(\vecq))$ satisfying (i)(ii) uniquely exists. With similar techniques, we show the existence of harmonic metrics for $SO(n,n+1)$-Higgs bundles in Collier's component and $Sp(4,\mathbb R)$-Higgs bundles in Gothen's component over $X$, under some mild assumptions.
    
\vspace{.1in}
\end{abstract}
\tableofcontents
\section{Introduction}
Let $X$ be a Riemann surface and $(E,\delbar_E,\theta)$ be a Higgs bundle on $X$.
Let $h$ be a Hermitian metric of $E$.
We obtain the Chern connection $\nabla_h=\delbar_E+\del_{E,h}$
and the adjoint $\theta^{*h}$ of $\theta$.
The metric $h$ is called a harmonic metric of
the Higgs bundle $(E,\delbar_E,\theta)$
if $\nabla_h+\theta+\theta^{*h}$ is flat,
i.e.,
$\nabla_h\circ\nabla_h+[\theta,\theta^{*h}]=0$.
It was introduced by Hitchin \cite{selfduality},
and it has been one of the most important
and interesting mathematical objects. 
A starting point is the study 
of the existence and the classification of harmonic metrics.
If $X$ is compact,
the results of Hitchin \cite{selfduality}
and Simpson \cite{s1} show that a Higgs bundle is polystable
of degree $0$ if and only if it admits a harmonic metric. Together with the work of Corlette \cite{Corlette} and Donaldson \cite{Donaldson}, one obtains the non-Abelian Hodge correspondence which says the moduli space of polystable $SL(n,\mathbb C)$-Higgs bundles is isomorphic to the representation variety of the surface group $\pi_1(S)$ into $SL(n,\mathbb C)$.
The study of harmonic metrics for Higgs bundles in the non-compact case was pioneered by Simpson \cite{s1,Simpson},
and pursued by Biquard-Boalch \cite{Biquard-Boalch}
and the second author \cite{Mochizuki-KH-Higgs}.

Let $\vecq=(q_2,\cdots, q_n)$, where $q_j$ is a holomorphic $j$-differential on $X$. One can naturally construct a Higgs bundle $(\hyperk_{X,n},\theta(\vecq))$ as follows. Let $K_X$ be the canonical line bundle of $X$.
The multiplication of $q_j$ induces the following morphisms:
\[
 K_X^{(n-2i+1)/2}\to
K_X^{(n-2i+2(j-1)+1)/2}\otimes K_X
 \quad
 (j\leq i\leq n).
\]
We also have the identity map
for $i=1,\ldots,n-1$:
\[
 K_X^{(n-2i+1)/2}
 \to
 K_X^{(n-2(i+1)+1)/2}\otimes K_X.
\]
They define a Higgs field
$\theta(\vecq)$ of $\hyperk_{X,n}=\oplus_{i=1}^nK_X^{(n+1-2i)/2}$.
The natural pairings
$K_X^{(n-2i+1)/2}
\otimes
K_X^{-(n-2i+1)/2}
\to \nbigo_X$
induce a non-degenerate symmetric bilinear form
$C_{\hyperk,X,n}$ of
$\hyperk_{X,n}$. There exists a basis of $SL(n,\mathbb C)$-invariant homogeneous polynomials $p_i$ of deg $i (i=2,\cdots, n)$ on $sl(n,\mathbb C)$ such that $p_i(\theta(\vecq))=q_i$. The Hitchin fibration is from the moduli space of polystable $SL(n,\mathbb C)$-Higgs bundles to the vector space $\oplus_{i=2}^nH^0(X,K_X^i)$ given by $$[(E,\theta)]\longmapsto (p_2(\theta),\cdots, p_n(\theta)).$$ Such Higgs bundles $(\hyperk_{X,n},\theta(\vecq))$ were introduced by Hitchin in \cite{Hit92} for compact hyperbolic Riemann surfaces. They form a section of the Hitchin fibration. For this reason, for arbitrary (not necessarily compact) Riemann surfaces, we call  $(\hyperk_{X,n},\theta(\vecq))$  Higgs bundles in the Hitchin section.

 For the compact hyperbolic surface case, Hitchin in \cite{Hit92} showed that $(\hyperk_{X,n},\theta(\vecq))$ are always stable and the Hitchin section corresponds to Hitchin component, a connected component in the representation variety of $\pi_1(X)$ into $SL(n,\mathbb R)$ which contains embedded Fuchsian representations. In particular, when $n=2$, the Hitchin section parametrize the Teichm\"uller space. Hitchin component has been the central object in the field of higher Teichm\"uller theory. For the case when $X=\bar X-D$ where $\bar X$ is a compact Riemann surface and $D$ is a finite set of points, let $q_j(j=2,\cdots, n)$ be meromorphic differentilas on $\bar X$ with possible poles at $D$ of pole order at most $j-1$. Using the work of Simpson \cite{SimpsonNoncompact} on parabolic Higgs bundles, Biswas-Ar\'{e}s-Gastesi-Govindarajan in \cite{BiswasParabolicTeichmuller} showed $(\hyperk_{X,n},\theta(\vecq))$ can be prolonged to a stable parabolic Higgs bundle of degree $0$ over $\bar X$ and thus admits a harmonic metric. Moreover, the Hitchin section corresponds to a connected component of the representation variety of $\pi_1(X)$ into $SL(n,\mathbb R)$ such that the holonomy of loops around punctures are of certain parabolic holonomy. \\ 

We want to study Higgs bundles in the Hitchin section in general case: tuples of holomorphic differentials on an arbitrary non-compact Riemann surfaces, e.g., unit disk, of infinite topology, etc. We focus on the following natural question. 
\begin{question}\label{Question1}
Given a tuple of holomorphic differentials $\vecq=(q_2,\cdots, q_n)$ on a non-compact Riemann surface $X$,\\
(1) does there exist a harmonic metric on $(\hyperk_{X,n}, \theta(\vecq))$ compatible with $C_{\hyperk,X,n}$?\\
(2) If so, can one find a notion of ``best" harmonic metric such that it uniquely exists?
\end{question}
\begin{rem}
1. When $X$ is parabolic, that is, $\mathbb C$ or $\mathbb C^*$, there exists no harmonic metric on $(\hyperk_{X,n}, \theta(\mathbf{0}))$. When $X$ is hyperbolic, each hyperbolic K\"ahler metric over $X$ induces a harmonic metric on $(\hyperk_{X,n}, \theta(\mathbf{0}))$.\\
2. Suppose $n=2$, $q_2\neq 0$ and $X$ is an arbitrary non-compact Riemann surface. The work of \cite{WAN} \cite{AuWan} \cite{Li17} together show there uniquely exists a harmonic metric $h$ of unit determinant of $(\hyperk_{X,n}, \theta(q_2)) $ satisfying $(h|_{K_X^{-1/2}})^2$ defines a complete metric on $X$.\\
3. Suppose $\vecq=(0,\cdots,0,q_n)$, $q_n\neq 0$ and $X$ is an arbitrary non-compact Riemann surface, the authors in \cite{LiMochizuki} introduce the notion of a complete metric $h$ on  $(\hyperk_{X,n}, \theta(\vecq))$, that is $h$ is diagonal of unit determinant and satisfies $(h|_{K_X^{(n+1-2i)/2}})^{-1}\otimes (h|_{K_X^{(n+1-2(i+1))/2}}) (i=1,\cdots,n-1)$ defines a complete metric on $X$. And we show the existence and uniqueness for a complete metric of $(\hyperk_{X,n}, \theta(\vecq))$. Sagman \cite{Sagman} later extends the existence of complete metric to the subcyclic case $(0,\cdots, 0,q_{n-1},0)$. For such two cases of lower ranks, there are rich related geometry including hyperbolic affine spheres in $\mathbb R^3$ \cite{Labourie, Loftin0}, maximal surfaces in $\mathbb H^{2,2}$ \cite{CTT}, $J$-complex curves in $\mathbb H^{4,2}$\cite{g2geometry,NieXin,HolomorphicCT}. There are extensive studies on the harmonic metrics for such two cases over non-compact surfaces, see e.g. \cite{BenoistHulinFiniteVolume,BenoistHulin,DumasWolf, nie2019poles, TamburelliWolf, Parker, GuestLin, GuestItsLin, MochizukiToda1, MochizukiToda2}.\\
4. In \cite{LiMochizukiGeneric}, the authors consider generically regular semisimple Higgs bundles which admit a non-degenerate symmetric pairing $C$. Here the condition ``generically regular semisimple" means there exists a point such that the Higgs field has $n$ distinct eigen $1$-forms. For such Higgs bundles, the authors show the existence of a harmonic metric compatible with $C$. Note that this result is not restricted to Higgs bundles in the Hitchin section. 
\end{rem}

 A harmonic metric on $(\hyperk_{X,n}, \theta(\vecq))$ compatible with $C_{\hyperk,X,n}$ gives rise to a representation $\rho:\pi_1(X)\rightarrow SL(n,\mathbb R)$ and a $\rho$-equivariant harmonic map to the symmetric space $SL(n,\mathbb R)/SO(n)$. Here $SL(n,\mathbb R)/SO(n)$ is equipped with the $SL(n,\mathbb R)$-invariant Riemannian metric induced by the Killing form $B(X,Y)=2n\cdot \tr(XY)$ on $sl(n,\mathbb R)$.  A closely related question is as follows. 
\begin{question}\label{Question2}
Given a tuple of holomorphic differentials $\vecq=(q_2,\cdots, q_n)$ on a non-compact Riemann surface $X$, does there exist an equivariant harmonic map $f:\widetilde X\rightarrow SL(n,\mathbb R)/SO(n)$ such that $q_i=p_i(-\frac{1}{2}f^{-1}\partial f)$ for $i=2,\cdots, n$? 
\end{question} Here we used the explicit relation $-\frac{1}{2}f^{-1}\partial f=\theta(\vecq),$ see e.g. \cite[Section 5.1]{LiIntroduction}. If Question \ref{Question1}(1) holds for some $\vecq_0$ on $X$, then Question \ref{Question2} automatically holds for $\vecq_0$ on $X$. 
\begin{rem}
When $X=\mathbb C$ and $\vecq$ are polynomial differentials, Question \ref{Question2} reduces to the question of Tamburelli-Wolf in \cite[Question A]{TamburelliWolf}.
\end{rem}

\subsection{Harmonic metrics for Higgs bundles in the Hitchin section}
Suppose $X$ is a non-compact hyperbolic Riemann surface, equivalently, it is not $\mathbb C$ nor $\mathbb C^*$ . Let $g_X$ be the unique complete hyperbolic K\"ahler metric on $X$. Let $h_X=\oplus_{k=1}^na_k\cdot g_X^{-\frac{n+1-2k}{2}}$, where $a_k$ are some fixed constants. Such $a_k$'s are chosen so that $h_X$ is a harmonic metric for the Higgs bundle $(\hyperk_{X,n}, \theta(\mathbf{0}))$. 

Let $F_k=\oplus_{l\leq k}K_X^{\frac{n+1-2l}{2}}.$ Then $\{0\subset F_1\subset F_2\subset\cdots \subset F_n\}$ forms an increasing filtration of $\hyperk_{X,n}$. We call a Hermitian metric $h$ on $\hyperk_{X,n}$ \textbf{weakly dominates} $h_X$ if $\det(h|_{F_k})\leq \det(h_X|_{F_k})$ for $1\leq k\leq n-1.$ Our main result in this paper is the following two theorems, as an answer to Question \ref{Question1}.
\begin{thm}(Theorem \ref{existence})
On a non-compact hyperbolic surface $X$, there exists a harmonic metric $h$ on $(\hyperk_{X,n},\theta(\vecq))$ satisfying (i) $h$ weakly dominates $h_X$; (ii) $h$ is compatible with $C_{\hyperk,X,n}.$

As a result, the associated harmonic map $f: (\widetilde X,\widetilde{g_X})\rightarrow SL(n,\mathbb R)/SO(n)$ satisfies the energy density $e(f)\geq \frac{n^2(n^2-1)}{6}.$ The equality holds if $\vecq=0.$
\end{thm}

\begin{thm} (Theorem \ref{BoundedDifferentialExistence})\label{UniquenessIntro}
On a non-compact hyperbolic surface $X$, suppose $q_i (i=2,\cdots,n)$ are bounded with respect to $g_X$. Then there uniquely exists a harmonic metric $h$ on $(\hyperk_{X,n},\theta(\vecq))$ satisfying (i) $h$ weakly dominates $h_X$; (ii) $h$ is compatible with $C_{\hyperk,X,n}.$

Moreover, $h$ is mutually bounded with $h_X.$
\end{thm}

As an application of Theorem \ref{UniquenessIntro}, we reprove the existence and uniqueness of a harmonic metric on $(\hyperk_{X,n},\theta(\vecq))$ over a compact hyperbolic Riemann surface. Note that our proof here does not invoke the Hitchin-Kobayashi correspondence by using the stability of Higgs bundle. 
\begin{thm}\label{CompactExistenceIntro} (Theorem \ref{CompactExistence})
Given a tuple of holomorphic differentials $\vecq=(q_2,\cdots, q_n)$ on a compact hyperbolic surface $X$, there uniquely exists a harmonic metric $h$ on $(\hyperk_{X,n},\theta(\vecq))$ satisfying $h$ is compatible with $C_{\hyperk,X,n}.$

Moreover, $h$ weakly dominates $h_X$.
\end{thm}

\subsection{Harmonic metrics for Higgs bundles which admit a full filtration}
In fact, we prove the existence of harmonic metrics for a more general family of Higgs bundles than Higgs bundles in the Hitchin section. Consider a Higgs bundle $(E,\theta)$ over a Riemann surface $X$ which admits a full holomorphic filtration $\mathbf F=\{0\subset F_1\subset F_2\subset \cdots\subset F_n\}$. We require that the induced map $\theta$ on each $Gr_k(E):=F_k/F_{k-1}$ is not a zero map, $k=1,\cdots, n-1$. Let $E_0=\oplus_{k=1}^nGr_k(E)$ and $\theta_0$ are induced by $\theta$ on the graded bundles $Gr_k(E)$. Then $(E_0, \theta_0)$ is a holomorphic chain of type $(1,1,\cdots,1).$ Take $F_k^0=\oplus_{l\leq k}Gr_l(E)$. 
There is a canonical way identifying $\det(F_k)$ and $\det(F_k^0)$, for $1\leq k\leq n$. So a metric on $\det(F_k)$ can be viewed as a metric on $\det(F_k^0)$. 

\begin{df}
Let $h,h_1$ be Hermitian metrics on $E, E_0$ respectively. Call $h$ weakly dominates $h_1$ if $$\det(h|_{F_k})\leq \det(h_1|_{F_k^0}),\quad 1\leq k\leq n-1.$$
\end{df}

We prove the following existence result.
\begin{thm}(Theorem \ref{MainTheorem})\label{MainTheoremIntro}
Suppose there exists a diagonal harmonic metric $h_1$ on $(E_0,\theta_0)$, then there exists a harmonic metric $h$ on $(E, \theta)$ satisfying (i) $\det(h)=\det(h_1)$; (ii) $h$ weakly dominates $h_1$.
\end{thm}

Because of Theorem \ref{MainTheoremIntro}, we are interested in the existence of a diagonal harmonic metric on a holomorphic chain of type $(1,\cdots,1).$ However, we find that such metric does not always exist, see Proposition \ref{boundedRank2} and Proposition \ref{SomeNoneExistence}. In Theorem \ref{NilpotentHiggsBundles}, we provide a sufficient condition of the existence of a harmonic metric on holomorphic chains. 

\subsection{$SO(n,n+1)$-Higgs bundles and $Sp(4,\mathbb R)$-Higgs bundles}
The Higgs bundles we consider in Theorem \ref{MainTheoremIntro} also appear in $SO(n,n+1)$-Higgs bundles in Collier section and $Sp(4,\mathbb R)$-Higgs bundles in Gothen section. As applications of Theorem \ref{MainTheoremIntro} and the existence result for diagonal harmonic metric on holomorphic chains, we show in \S\ref{SO(n,n+1)} the existence of harmonic metric on $SO(n,n+1)$-Higgs bundles in Collier section. In \S\ref{Sp(4,R)}, we show the existence of harmonic metric on $Sp(4,\mathbb R)$-Higgs bundles in Gothen section.

\subsection{Further questions}
~~~~ 1. Our techniques here only apply to hyperbolic Riemann surfaces since it relies on the existence of harmonic metric on the graded Higgs bundle. Since the graded Higgs bundles are nilpotent, the existence of a harmonic metric forces the Riemann surface to be hyperbolic. Therefore, $\vecq\neq \mathbf{0}$ is a necessary condition for the existence of harmonic metric on $(\hyperk_{X,n},\theta(\vecq))$ over a parabolic Riemann surface. So it would be interesting to ask if $\vecq\neq \mathbf{0}$ is a sufficient condition. So far, the best answer we can provide are Higgs bundles satisfying generically regular semisimple condition.

\vspace{0.15cm}
2. We would like to see the uniqueness result in Theorem \ref{UniquenessIntro} extends to all $\vecq$ without the boundedness condition.

\vspace{0.15cm}
3. For holomorphic chains, we find a sufficient condition for the existence of a harmonic metric in Theorem \ref{NilpotentHiggsBundles}. We would like to find a sufficient and necessary condition for the existence of a diagonal harmonic metric for holomorphic chains of type $(1,\cdots, 1)$. 

\vspace{0.15cm}
4. There is a natural $\mathbb C^*$-action on the space of gauge equivalent classes of Higgs bundles as follows: $t\cdot [(E,\theta)]=[(E,t\theta)].$ We want to ask if the $\mathbb C^*$-action preserve the property admitting a harmonic metric. More precisely, suppose a Higgs bundle $(E,\theta)$ admits a harmonic metric, does there exist a harmonic metric on $(E,t\cdot \theta)$ for $t\in \mathbb C^*$? This is true if the base Riemann surface is compact hyperbolic since the stability is preserved by the $\mathbb C^*$-action. For non-compact Riemann surfaces, the answer is unclear. The evidence for this conjecture is that the properties in the two cases we can prove the existence of harmonic metrics are preserved by the $\mathbb C^*$-action: (1) Higgs bundles being in the Hitchin section; (2) generically regular semisimple and admits a non-degenerate symmetric pairing.

\subsection*{Organization}
In \S\ref{Preliminaries},
we give some results on the existence of harmonic metric
using exhaustion family of harmonic metrics of Dirichlet problem.
In \S\ref{ExistenceSection},
we study the existence of harmonic metric for the Higgs bundles
which admit a full holomorphic filtration.
In \S\ref{UniquenessSection},
we study the uniqueness of real harmonic metrics
of some Higgs bundles
which are mutually bounded with a canonically constructed metric.
In \S \ref{HitchinSection}, we apply the existence result
to Higgs bundles in the Hitchin section and
show the uniqueness result for the case of bounded differentials.
In \S\ref{BoundedExistenceSection}, we show the existence of harmonic metric under boundedness condition on the Higgs bundle and apply it to holomorphic chains. In the last two sections, we show the existence of harmonic metric on $SO(n,n+1)$-Higgs bundles and  $Sp(4,\mathbb R)$-Higgs bundles.

\subsection*{Acknowledgement}
The first author is partially supported by the National Key R\&D Program of China No. 2022YFA1006600, the Fundamental Research
Funds for the Central Universities and Nankai Zhide foundation.
The second author is partially supported by
the Grant-in-Aid for Scientific Research (A) (No. 21H04429),
the Grant-in-Aid for Scientific Research (A) (No. 22H00094),
the Grant-in-Aid for Scientific Research (A) (No. 23H00083),
and the Grant-in-Aid for Scientific Research (C) (No. 20K03609),
Japan Society for the Promotion of Science.
He is also partially supported by the Research Institute for Mathematical
Sciences, an International Joint Usage/Research Center located in Kyoto
University.

\section{Preliminaries on existence of harmonic metrics}\label{Preliminaries}
In this section, we give some results on the existence of harmonic metric using exhaustion family of harmonic metrics of Dirichlet problem. A variant version also appears in \cite[Section 2]{LiMochizuki1}.

\subsection{Dirichlet problem}

Let $X$ be any Riemann surface.
Let $(E,\delbar_E,\theta)$ be a Higgs bundle on $X$.
For a Hermitian metric $h$ of $E$,
we obtain the Chern connection
$\nabla_h=\delbar_E+\del_E^h$
of $(E,\delbar_E,h)$.
The curvature of $\nabla_h$ is denoted by $F(\nabla_h)$ or $F(h)$.
We also obtain the adjoint
$\theta^{*h}$ of $\theta$
with respect to $h$.
The curvature of
$\nabla_h+\theta+\theta^{*h}$
is denoted by $F(E,\delbar_E,\theta)$,
i.e.,
$F(E,\delbar_E,\theta)=F(\nabla_h)+[\theta,\theta^{*h}]$.

Let $Y\subset X$ be a relatively compact connected open subset
with smooth boundary $\del Y$.
Assume that $\del Y$ is non-empty.
Let $h_{\del Y}$ be any Hermitian metric of
$E_{|\del Y}$.

\begin{prop}[Donaldson]
\label{prop;20.5.29.20}
 There exists a unique harmonic metric $h$
 of $(E,\delbar_E,\theta)$
 such that $h_{|\del Y}=h_{\del Y}$.
\end{prop}
\pf
This was proved by Donaldson
\cite[Theorem 2]{Donaldson-boundary-value}
in the case $Y$ is a disc.
The general case is essentially the same.
We explain an outline of the proof
for the convenience of the reader.

We may assume that $X$ is an open Riemann surface.
According to \cite{Gunning-Narasimhan},
there exists a nowhere vanishing
holomorphic $1$-form $\tau$ on $X$.
Let $f$ be the automorphism of $E$
determined by $\theta=f\,\tau$.
We consider the K\"ahler metric $g_X=\tau\,\taubar$ of $X$.

 Let $\Gamma$ be a lattice of $\cnum$
 and let $T$ be a real $2$-dimensional torus
 obtained as $\cnum/\Gamma$.
 We set $g_T=dz\,d\zbar$.
 We set $\Xtilde:=X\times T$
 with the projection $p:\Xtilde\lrarr X$.
 It is equipped with the flat K\"ahler metric
 $g_{\Xtilde}$ induced by $g_T$ and $g_X$.
 We set $\Ytilde:=p^{-1}(Y)$.

 Let $\Etilde$ be the pull back of $E$
 with the holomorphic structure
 $p^{\ast}(\delbar_{E})+p^{\ast}(f)\,d\zbar$.
According to the dimensional reduction of Hitchin,
a Hermitian metric $h$ of $E_{|Y}$ is
a harmonic metric of $(E,\delbar_E,\theta)_{|Y}$
if and only if
$\Lambda_{\Ytilde}F(p^{\ast}h)=0$.
According to a theorem of
Donaldson \cite[Theorem 1]{Donaldson-boundary-value},
there exists a unique Hermitian metric $\htilde$
of $\Etilde$
such that
$\Lambda_{\Ytilde}F(\htilde)=0$
and that
$\htilde_{|p^{-1}(\del Y)}=p^{\ast}(h_{\del Y})$.
By the uniqueness,
$\htilde$ is $T$-invariant.
Hence, there uniquely exists a harmonic metric $h$
of $(E,\delbar_{E},\theta)_{|Y}$
which induces $\htilde$.
It satisfies
$h_{|\del Y}=h_{\del Y}$.
\hfill\qed
\vspace{.1in}

Let $h_0$ be a Hermitian metric of $E$.
Assume that $\det(h_0)$ is flat.
 \begin{cor}
There exists a unique harmonic metric $h$ of
  $E_{|Y}$ such that
  $h_{|\del Y}=h_{0|\del Y}$
  and that
  $\det(h)=\det(h_0)_{|Y}$.
 \end{cor}
 \pf
 There exists a unique harmonic metric $h$
 such that $h_{|\del Y}=h_{0|\del Y}$.
 We obtain
 $\det(h)_{|\del Y}=\det(h_0)_{|\del Y}$.
 Note that 
 both $\det(h)$ and $\det(h_0)_{|Y}$ are flat.
 By the uniqueness in Proposition \ref{prop;20.5.29.20},
 we obtain $\det(h)=\det(h_0)_{|Y}$. 
 \hfill\qed

\subsection{Convergence}
\label{subsection;23.6.16.10}

Let $X$ be an open Riemann surface.
Let $h_0$ be a Hermitian metric of $E$.

\begin{df}
 An exhaustive family  $\{X_i\}$ of a Riemann surface $X$
 means an increasing sequence of relatively compact open subsets
 $X_1\subset X_2\subset\cdots$ of $X$
 such that $X=\bigcup X_i$.
 The family is called smooth if
 $\del X_i$ are smooth.
 \hfill\qed
\end{df}

Let $\{X_i\}$ be a smooth exhaustive family of $X$.
The restriction $h_{0|X_i}$ is denoted by $h_{0,i}$.
Let $h_i$ $(i=1,2,\ldots)$
be harmonic metrics of
$(E,\delbar_{E},\theta)_{|X_i}$.
Let $s_i$ be the automorphism of
$E_{|X_i}$ determined by
$h_i=h_{0,i}\cdot s_i$.
Let $f$ be an $\real_{>0}$-valued function on $X$
such that each $f_{|X_i}$ is bounded. Though the following proposition is proved in \cite{LiMochizuki1}, we include the proof for the
convenience of the readers.
\begin{prop}
\label{prop;20.5.29.1}
Assume that 
 $|s_i|_{h_{0,i}}+|s^{-1}_i|_{h_{0,i}}\leq f_{|X_i}$
 for any $i$.
 Then, there exists a subsequence
 $s_{i(j)}$ which is convergent
 to an automorphism $s_{\infty}$ of $E$
 on any relatively compact subset of $X$
 in the $C^{\infty}$-sense.
 As a result, we obtain a harmonic metric
 $h_{\infty}=h_0s_{\infty}$ of $(E,\delbar_E,\theta)$
 as the limit of the subsequence $h_{i(j)}$.
 Moreover, we obtain
 $|s_{\infty}|_{h_0}+|s_{\infty}^{-1}|_{h_0}\leq f$.
 In particular, if $f$ is bounded, $h_0$ and $h_{\infty}$
 are mutually bounded.
\end{prop}
\pf
We explain an outline of the proof.
Let $g_X$ be a K\"ahler metric of $X$.
According to a general formula (\ref{eq;20.8.16.3}) below,
the following holds on any $X_i$:
\begin{equation}
\label{eq;20.6.30.1}
\sqrt{-1}\Lambda\delbar\del\tr(s_i)
=
-\sqrt{-1}\tr\bigl(
s_i\Lambda F(h_{0,i})
\bigr)
-\bigl|(\delbar+\theta)(s_i)\cdot s_i^{-1/2}\bigr|^2_{h_{0,i},g}.
\end{equation}

Let $K$ be any compact subset of $X$.
Let $N$ be a relatively compact neighbourhood of $K$ in $X$.
Let $\chi:X\lrarr\real_{\geq 0}$ be a $C^{\infty}$-function
such that (i) $\chi_{|K}=1$,
(ii) $\chi_{|X\setminus N}=0$,
(iii) $\chi^{-1/2}\del\chi$ and $\chi^{-1/2}\delbar\chi$
on $\{P\in X\,|\,\chi(P)>0\}$ induces a $C^{\infty}$-function on $X$.

There exist $i_0$ such that
$N$ is a relatively compact open subset of $X_i$
for any $i\geq i_0$.
We obtain the following:
\begin{equation}
\label{eq;20.6.30.2}
 \sqrt{-1}\Lambda\delbar\del(\chi\tr(s_i))
 =\chi\sqrt{-1}\Lambda\delbar\del\tr(s_i)
 +(\sqrt{-1}\Lambda\delbar\del\chi)\cdot\tr(s_i)
 +\sqrt{-1}\Lambda(\delbar\chi\del\tr(s_i))
 -\sqrt{-1}\Lambda(\del\chi\delbar\tr(s_i)).
\end{equation}
Note that
$|\delbar_Es_i|_{h_i,g_X}
=|\del_{E,h_i}s_i|_{h_i,g_X}$,
and that
\begin{equation}
\label{eq;20.6.30.3}
\left|
 \int_X\sqrt{-1}\Lambda(\delbar\chi\del\tr(s_i))
 \right|
 \leq
 \left(
 \int_X|\chi^{-1/2}\delbar\chi|^2
 \right)^{1/2}
 \cdot
 \left(
 \int_X \chi|\del_{E,h_i} s_i|^2_{h_{0},g_X}
 \right)^{1/2}.
\end{equation}
Note that there exists $C_0>0$
such that $|s_i|_{h_0}+|s_i^{-1}|_{h_0}\leq C_0$ on $N$.
By (\ref{eq;20.6.30.1}), (\ref{eq;20.6.30.2})
and (\ref{eq;20.6.30.3}),
there exist $C_j>0$ $(j=1,2)$ such that
the following holds for any sufficiently large $i$:
\[
 \int\chi
 \bigl|
 \delbar_Es_i
 \bigr|^2_{h_{0},g_X}
 +\int\chi\bigl|[\theta, s_i]\bigr|^2_{h_{0},g_X}
 \leq
 C_1+
 C_2\left(
 \int
 \chi\bigl|
 \delbar_Es_i
 \bigr|^2_{h_{0},g_X}
 +\int\chi\bigl|[\theta, s_i]\bigr|^2_{h_{0},g_X}
 \right)^{1/2}
\]
Therefore, there exists $C_3>0$ such that
the following holds for any sufficiently large $i$:
\[
\int\chi
 \bigl|
 \delbar_Es_i
 \bigr|^2_{h_{0},g_X}
 +\int\chi\bigl|[\theta, s_i]\bigr|^2_{h_{0},g_X}
 \leq C_3.
\]
We obtain the boundedness of the $L^2$-norms of
$\delbar_Es_i$
and $\del_{E,h_i}s_i$ $(i\geq i_0)$
on $K$
with respect to $h_0$ and $g_X$.
By a variant of Simpson's main estimate
(see \cite[Proposition 2.1]{MochizukiAsymptotic}),
we obtain the boundedness of the sup norms of
$\theta$ on $N$ with respect to
$h_i$ and $g_X$.
By the Hitchin equation,
we obtain the boundedness of
the sup norms of 
$\delbar_E(s_i^{-1}\del_{E,h_i} s_i)$ on $N$
with respect to $h_i$ and $g_X$.
By using the elliptic regularity,
we obtain that
the $L_1^p$-norms of
$s_i^{-1}\del_{E,h_i}(s_i)$ on
a relatively compact neighbourhood of $K$
are bounded
for any $p>1$.
It follows that
$L_2^p$-norms of $s_i$ on
a relatively compact neighbourhood of $K$
are bounded for any $p$.
Hence,
a subsequence of $s_i$
is weakly convergent in $L_2^p$
on a relatively compact neighbourhood of $K$.
By the bootstrapping argument
using a general formula (\ref{eq;20.8.16.5}) below,
we obtain that the sequence is convergent
on a relatively compact neighbourhood of $K$
in the $C^{\infty}$-sense.
By using the diagonal argument,
we obtain that 
a subsequence of $s_i$ is weakly convergent
in $C^{\infty}$-sense on any compact subset.
\hfill\qed

\vspace{.1in}

\subsection{Appendix}
\label{subsection;20.8.16.10}

We recall some fundamental formulas due to
Simpson \cite[Lemma 3.1]{s1}
for the convenience of the readers.

Let $h_i$ $(i=1,2)$ be Hermitian metrics of $E$.
We obtain the automorphism $s$ of $E$ determined by
$h_2=h_1\cdot s$.
Let $g$ be a K\"ahler metric of $X$,
let $\Lambda$ denote the adjoint of the multiplication of
the associated K\"ahler form.
Then, according to \cite[Lemma 3.1 (a)]{s1},
we obtain the following on $X$:
\begin{equation}
\label{eq;20.8.16.5}
 \sqrt{-1}\Lambda
 \bigl(\delbar_E+\theta\bigr)
 \circ
 \bigl(\del_{E,h_1}
 +\theta^{*h_1}\bigr)s
=
 s\sqrt{-1}\Lambda\bigl(
 F(h_2)-F(h_1)
 \bigr)
+\sqrt{-1}\Lambda
 \Bigl(
 \bigl(\delbar_{E}+\theta\bigr)(s)
  s^{-1}
  \bigl(\del_{E,h_1}
  +\theta^{*h_1}\bigr)(s)
 \Bigr).
\end{equation}
By taking the trace,
and by using \cite[Lemma 3.1 (b)]{s1},
we obtain
\begin{equation}
\label{eq;20.8.16.3}
 \sqrt{-1}\Lambda\delbar\del
 \tr(s)
=
 \sqrt{-1}
 \tr\Bigl(
  s\Lambda\bigl(F(h_2)-F(h_1)\bigr)
 \Bigr)
-\Bigl|
 \bigl(\delbar_{E}+\theta\bigr)(s)
  s^{-1/2}
 \Bigr|^2_{h_1,g}.
\end{equation}
Note that
$(\delbar_{E}+\theta)(s)
=\delbar_{E}(s)
+[\theta,s]$.
Moreover,
$\delbar_E(s)$
is a $(0,1)$-form,
and
$[\theta,s]$
is a $(1,0)$-form.
Hence, (\ref{eq;20.8.16.3}) is also rewritten as follows:
\begin{equation}
\label{eq;20.8.16.2}
 \sqrt{-1}\Lambda\delbar\del\tr(s)
 =
 \sqrt{-1}
 \tr\Bigl(
  s\Lambda\bigl(F(h_2)-F(h_1)\bigr)
 \Bigr)
-\bigl|
[\theta,s]s^{-1/2}
\bigr|^2_{h_1,g}
-\bigl|
\delbar_E(s)s^{-1/2}
\bigr|_{h_1,g}.
\end{equation}
We also recall
the following inequality \cite[Lemma 3.1 (d)]{s1}:
\begin{equation}
\label{eq;20.8.16.6}
 \sqrt{-1}\Lambda\delbar\del\log\tr(s)
  \leq
  \bigl|\Lambda F(h_1)\bigr|_{h_1}
 +\bigl|\Lambda F(h_2)\bigr|_{h_2}.
\end{equation}
In particular,
if both $h_i$ are harmonic,
the functions $\tr(s)$
and $\log\tr(s)$ are subharmonic:
\begin{equation}
\label{eq;20.8.17.1}
 \sqrt{-1}\Lambda\delbar\del\tr(s)
  =-\bigl|
  (\delbar_E+\theta)(s)s^{-1/2}
  \bigr|^2_{h,g}\leq 0,
  \quad\quad
  \sqrt{-1} \Lambda\delbar\del\log\tr(s)
  \leq 0.
\end{equation}

\section{Domination property and the existence of harmonic metrics}\label{ExistenceSection}

\subsection{Full flags and Hermitian metrics}\label{FiltrationVectorSpace}
Let $V$ be a complex vector space equipped with a base $\mathbf{e}=(e_1,\cdots, e_n)$. For $k=1,\cdots, n$, let $F_k(V)$ denote the subspace generated by $e_1,\cdots, e_k.$ We set $F_0(V)=0$. We set $Gr_k^F(V)=F_k(V)/F_{k-1}(V)$. There exists a natural isomorphism 
$$\rho_k:Gr_k^F(V)\otimes \det(F_{k-1})\cong \det(F_k).$$
Let $h$ be a Hermitian metric of $V$. Let $F_k(h)$ denote the induced metric of $F_k(V).$ It induces a Hermitian metric $\det(F_k(h))$ of $\det(F_k(V)).$ Let $G_k(V,h)$ be the orthogonal complement of $F_{k-1}(V)$ in $F_k(V).$ The projection $F_k(V)\rightarrow Gr_k^F(V)$ induces an isomorphism $G_k(V,h)\cong Gr_k^F(V).$ We obtain the metric $Gr_k^F(h)$ of $Gr_k^F(V)$ which is induced by $h|_{Gr_k(V,h)}$ and the isomorphism $G_k(V,h)\cong Gr_k^F(V).$

\begin{lem}
$\rho_k$ is isometric with respect to $\det(F_k(h))$ and $Gr_k^F(h)\otimes \det(F_{k-1}(h)).$
\end{lem}
\pf
There exists the orthogonal decomposition 
$F_k(V)=\oplus_{j=1}^kG_j(V,h).$ We choose $v_j\in G_j(V,h)$ such that $h(v_j,v_j)=1.$ The norm of $v_1\wedge\cdots\wedge v_k$ with respect to $\det(F_k(h))$ is $1$. Let $[v_k]\in Gr_k^F(V)$ denote the element induced by $v_k$. The norm of $[v_k]$ with respect to $Gr_k^F(h)$ is $1$. Then we obtain the claim of the lemma.
\hfill\qed\\

Denote $F_k^0(V)=\oplus_{l\leq k}Gr_k^F(V).$ It has the induced metric $F_k^0(h)$ from $h$ on $V$. From $\rho_k$'s, one naturally has an isomorphism between $F_k(V)$ with $F_k^0(V)$, which is an isometry with respect to $\det(F_k(h))$ and $\det(F_k^0(h)).$

\subsection{Set-up}\label{Setup}
Let $X$ be a hyperbolic Riemann surface and $K$ be its canonical line bundle.

Consider a Higgs bundle $(E, \theta)$ over $X$
which admits a full holomorphic filtration
$$\mathbf F=\{0=F_0\subset F_1\subset F_2\subset \cdots\subset F_n=E\}$$
and $\theta:F_k\rightarrow F_{k+1}\otimes K$.
We require that the induced map $\theta$ on each $F_k/F_{k-1}$ is not a zero map,
denoted by $\phi_k$, for $1\leq k\leq n-1$.

Denote by $Gr_k^F(E)$ the quotient line bundle $F_k/F_{k-1}$, equipped with the quotient holomorphic structure. Consider the holomorphic vector bundle $E_0=Gr^F(E):=\oplus_{k=1}^nGr_k^F(E)$. Let $\theta_0$ be formed by $\phi_k:Gr_k^F(E)\rightarrow Gr_{k+1}^F(E)\otimes K,$ for $1\leq k\leq n-1$. Therefore, $(E_0, \theta_0)$ is a holomorphic chain of type $(1,1,\cdots,1).$ Let $F_k^{(0)}=\oplus_{l\leq k}Gr_l^F(E).$ 

Let $h$ be a Hermitian metric on $E$. Let $F_k(h)$ denote the induced metric of $h$ on $F_k$. The metric $h$ induces a metric $Gr_k^F(h)$ on each $Gr_k^F(E)$, a diagonal metric $F_k^0(h)$ on $F_k^0$, and a diagonal metric on $E_0$. 
\begin{df}
Suppose $h$ is a Hermitian metric on $E$,
 and $h_1$ is a diagonal Hermitian metrics on $E_0$.
 Call $h$ weakly dominates $h_1$ if 
 \begin{equation}\label{Condition1}
   \det(F_k^0(h))\leq \det(F_k^0(h_1)),\quad 1\leq k\leq n-1.
 \end{equation}
\end{df}
Under the natural identification between $\det(F_k^0)$ and $\det(F_k)$ in \S\ref{FiltrationVectorSpace}, we can write the condition (\ref{Condition1}) as follows 
\begin{equation}\label{Condition2}
        \det(F_k(h))\leq \det(F_k^0(h_1)),\quad 1\leq k\leq n-1.
\end{equation}

\subsubsection{The graded case}
\label{subsection;23.3.31.2}

If $E=\oplus_{i=1}^nL_i$ and $F_k=\oplus_{l\leq k}L_l$
for some holomorphic line bundles $L_i$ over $X$,
there is a canonical isomorphism between $E$ and $E_0$
by mapping $L_k$ to $Gr_k^F(E)$.
In this case, we can identify $E$ and $E_0$
and view the metric $h_1$ as a metric on $E$ too.
Then we call $h$ weakly dominates $h_1$ if 
\begin{equation}\label{Condition3}
\det(F_k(h))\leq \det(F_k(h_1)),\quad 1\leq k\leq n-1.
 \end{equation}

Because of the following lemma,
we may assume the existence of such a grading
if the Riemann surface is non-compact.

\begin{lem}
Let $E$ be a holomorphic vector bundle  of rank $r$
on a non-compact Riemann surface $X$
equipped with an increasing filtration $F_j(j=1,\cdots, n)$
such that $\rank \Gr^F_j(E)=1$ for $j=1,\cdots, n$.
Then, there exists a frame $e_1,\cdots, e_n$ of $E$
such that $F_j=\oplus_{i\leq j}\mathcal O_Xe_i$.
\end{lem}
\pf
It is well known that $H^1(X,\mathcal O_X)=0$.
Because any holomorphic vector bundle $E'$ on $X$ is
isomorphic to $\mathcal O_X^{\rank(E')}$,
we have $H^1(X,E')=0$.
Let $\pi:E_1\rightarrow E_2$ be an epimorphism of holomorphic vector bundles on $X$.
Let $K$ be the kernel. Because 
\[H^0(X, \Hom(E_2,E_1))\rightarrow H^0(X, \Hom(E_2,E_2))\rightarrow H^1(X,\Hom(E_2,K))=0\] is exact, there exists a splitting
$s:E_2\rightarrow E_1$ such that $\pi\circ s=id_{E_2}$.
Then, the claim of the lemma follows.
\hfill\qed

\subsubsection{Symmetric pairings}
Let us recall the notion of compatibility of
a non-degenerate symmetric pairing
and a Hermitian metric
on a complex vector space $V$.
(See \cite[\S2.1]{LiMochizukiGeneric} for more details.)
Let $V^{\lor}$ denote the dual space of $V$.
Let $\langle\cdot,\cdot\rangle:V^{\lor}\times V\to \cnum$
denote the canonical pairing.

Let $C:V\times V\to\cnum$
be a non-degenerate symmetric bilinear form.
We obtain the linear isomorphism
$\Psi_C:V\simeq V^{\lor}$
by
$\langle \Psi_C(u),v\rangle=C(u,v)$.
We obtain the symmetric bilinear form
$C^{\lor}:V^{\lor}\times V^{\lor}\to\cnum$
by
\[
 C^{\lor}(u^{\lor},v^{\lor})
 =C(\Psi_C^{-1}(u^{\lor}),\Psi_C^{-1}(v^{\lor})).
\]
We have $\Psi_{C^{\lor}}\circ\Psi_C=\id_V$.

Let $h$ be a Hermitian metric of $V$.
We obtain the sesqui-linear isomorphism
$\Psi_h:V\simeq V^{\lor}$
by
$\langle\Psi_h(u),v\rangle=h(v,u)$.
We obtain the Hermitian metric $h^{\lor}$ of $V^{\lor}$
by
\[
h^{\lor}(u^{\lor},v^{\lor})=
h\bigl(
 \Psi_h^{-1}(v^{\lor}),\Psi_h^{-1}(u^{\lor})
\bigr).
\]
It is easy to see that
$\Psi_{h^{\lor}}\circ\Psi_h=\id_V$.
\begin{df}
We say that $h$ is compatible with $C$
if $\Psi_C$ is isometric with respect to $h$ and $h^{\lor}$.
 \hfill\qed
\end{df}
\begin{lem}
\label{lem;22.8.29.10}
The following conditions are equivalent.
\begin{itemize}
 \item  $h$ is compatible with $C$.
 \item $C(u,v)=\overline{C^{\lor}(\Psi_h(u),\Psi_h(v))}$ holds
       for any $u,v\in V$.
 \item $\Psi_{C^{\lor}}\circ\Psi_h=\Psi_{h^{\lor}}\circ\Psi_C$ holds.
It is also equivalent to
$\Psi_{C}\circ\Psi_{h^{\lor}}
=\Psi_h\circ\Psi_{C^{\lor}}$.
\end{itemize}
\end{lem}

Note that $h$ and $C$ induce a Hermitian metric $\det(h)$ and a non-degenerate symmetric pairing $\det(C)$ of $\wedge^nV$ respectively. The following lemma is clear. 
\begin{lem}
If $h$ is compatible with $C$, then $|\det(h)|=|\det(C)|$.
\end{lem}

\subsubsection{Symmetric pairing and graded bundles}\label{SymmetricForms}
Consider the graded case $E=\oplus_{i=1}^nL_i$ and $F_k=\oplus_{l\leq k}L_l$
for some holomorphic line bundles $L_i$ over $X$. Suppose in addition $L_i=L_{n+1-i}^{-1}$. 
Then it induces a natural symmetric pairing induced by $L_i\otimes L_{n+1-i}\rightarrow \mathcal O$, denoted by $C$.
In this case $\det(E)\cong \mathcal O$ and $|\det(C)|=1$. 
If $h$ is compatible with $C$, then $\det(h)=1$. 

With respect to the decomposition $\End(E)=\oplus_{ij} \Hom(L_i,L_j)$, the Higgs field $\theta=\sum _{ij}\theta_{ij}$ . 

Suppose moreover $\theta$ is symmetric with respect to $C$. That is, $\theta_{ij}=\theta_{n+1-j,n+1-i}$ under the identification between $\Hom(L_i,L_j)\cong \Hom (L_{n+1-j}, L_{n+1-i})=\Hom(L_j^{-1},L_i^{-1})$.

The graded bundle $E_0$ has an induced pairing that $C_0(Gr_k^FE, Gr_l^FE)=\delta_{kl}$. The canonical isomorphism between $E$ and $E_0$ takes $C$ to $C_0$. So we identify $(E_0,C_0)$ with $(E,C)$. Since $\theta_0=\sum_{i=1}^{n-1}\theta_{i,i+1}$, it is again symmetric with respect to $C$.

\subsection{Domination property and the Dirichlet problem}
Let $X$, $(E,\theta)$ and $(E_0,\theta_0)$ be as in \S\ref{Setup}. Let $h_1$ be a harmonic metric on $(E_0,\theta_0)$ orthogonal to the decomposition $E_0=\oplus_{k=1}^nGr_k^F(E)$. The following proposition is motivated by the result for Higgs bundles in the Hitchin section over compact Riemann surfaces in \cite{Li19}. Here we extend the result to surfaces with boundaries and more general Higgs bundles. This domination property turns out to be the key property in showing the convergence of harmonic metrics in the exhaustion process. 
\begin{prop}\label{Control} 
On a Riemann surface $X$ with boundary $\partial X$, suppose $(E,\theta)$ has a harmonic metric $h$ satisfying $h=h_1$ on $\partial X.$ Then $h$ weakly dominates $h_1.$
\end{prop}
\pf
For a holomorphic subbundle $F$ of $E$, we would like to deduce the Hitchin equation which respects $F$. Denote by $F^{\perp}$ the subbundle of $E$ perpendicular to $F$ with respect to the harmonic metric $H$. $F^{\perp}$ can be equipped with the quotient holomorphic structure from $E/F$. With respect to the $C^{\infty}$ orthogonal decomposition $$E=F\oplus F^{\perp},$$ we have the expression of the holomorphic structure $\delbar_E$ and the Higgs field $\phi$ as follows:
\begin{equation*}
\delbar_E=\begin{pmatrix}\bar\partial_F&\beta\\
0&\bar\partial_{F^{\perp}}\end{pmatrix},\quad \theta=\begin{pmatrix}\phi_1&\alpha\\B&\phi_2\end{pmatrix},\quad H=\begin{pmatrix}H_1&0\\0&H_2\end{pmatrix}, 
\end{equation*}
where the term $B\in \Omega^{1,0}(X,\Hom(F, F^{\perp}))$, $\alpha\in \Omega^{1,0}(X, \Hom(F^{\perp},F))$, and $\beta\in \Omega^{0,1}(X, \Hom(F^{\perp},F))$.

The Chern connection $\nabla_H$ and the adjoint $\theta^{*_H}$ of the Higgs field are
\begin{equation*}
\nabla_H=\begin{pmatrix}
\nabla_{H_1}&\beta\\
-\beta^{*_H}&\nabla_{H_2}\end{pmatrix},\quad 
\theta^{*_H}=\begin{pmatrix}
\phi_1^{*_{H_1}}&B^{*_H}\\
\alpha^{*_H}&\phi_2^{*_{H_2}}\end{pmatrix}.
\end{equation*}

We calculate the Hitchin equation with respect to the decomposition $E=F\oplus F^{\perp}$ and by restricting to $\Hom(F, F)$, we obtain
\begin{equation*}
F(\nabla_{H_1})-\beta\wedge\beta^{*_H}+\alpha\wedge\alpha^{*_H}+B^{*_H}\wedge B+[\phi_1,\phi_1^{*_{H_1}}]=0.
\end{equation*}
By taking trace and noting that $\tr([\phi_1,\phi_1^{*_{H_1}}])=0$, we obtain
\begin{equation*}
\tr(F(\nabla_{H_1}))-\tr(\beta\wedge\beta^{*_H})+\tr(\alpha\wedge\alpha^{*_H})+\tr(B^{*_H}\wedge B)=0.
\end{equation*}

Let $g_X=g(z)(dx^2+dy^2)$ be a conformal Riemannian metric on $X$.  The associated K\"ahler form associated to $g_X$ is \[\omega=\frac{\sqrt{-1}}{2}g(z)dz\wedge d\bar z.\] Note that $$|\partial/\partial_z|_{g_X}^2=\frac{g(z)}{2}, \quad |dz|_{g_X}^2=\frac{2}{g(z)}.$$ Thus the induced Hermitian metric on $K_X^{-1}$ can be written as $\frac{g(z)}{2} dz\otimes d\bar z$, still denoted as $g_X$. Denote by $\Lambda_{g_X}$ the contraction with respect to the K\"ahler form $\omega$. Therefore, 
\begin{eqnarray}\label{standard}
-\sqrt{-1}\Lambda_{g_X} \tr(F(\nabla_{H_1}))-\sqrt{-1}\Lambda_{g_X} \tr(B^{*_H}\wedge B)&=&-\sqrt{-1}\Lambda_{g_X}\tr(\beta\wedge\beta^{*_H})+\sqrt{-1}\Lambda_{g_X} \tr(\alpha\wedge\alpha^{*_H})\nonumber\\
&=&||\beta||_{H,g_X}^2+||\alpha||_{H,g_X}^2\geq 0.
\end{eqnarray}

We will apply the above procedure to $F=F_k$ for each $k=1,2,\cdots,n-1$. We take $L_i$ to be the perpendicular line bundle of $F_{i-1}$ inside $F_i$ with respect to he harmonic metric $h.$
Then we have a smooth decomposition of $$E=L_1\oplus L_2\oplus \cdots\oplus L_n.$$ With respect to the decomposition, we have the following:

I. the Hermitian metric $H$ solving the Hitchin equation is given by 
\begin{equation}\label{Metric}
H=\begin{pmatrix}h|_{L_1}&&&\\&h|_{L_2}&&\\&&\ddots&\\&&&h|_{L_n}\end{pmatrix}\end{equation}
where $h|_{L_i}$ is the induced Hermitian metric on $L_i$ and $h|_{L_i}=\det(h|_{F_i})/\det(h|_{F_{i-1}})$;

II. the holomorphic structure on $E$ is given by the $\bar\partial$-operator \begin{eqnarray}
\delbar_E=\begin{pmatrix}\bar\partial_1&\beta_{12}&\beta_{13}&\cdots&\beta_{1n}\\&\bar\partial_2&\beta_{23}&\cdots&\beta_{2n}\\&&\bar\partial_3&\cdots&\beta_{3n}\\&&&\ddots&\vdots\\&&&&\bar\partial_n\end{pmatrix}\label{HolomorphicStructure}
\end{eqnarray} where $\bar\partial_k$ are $\bar\partial$-operators defining the holomorphic structures on $L_k$, and $\beta_{ij}\in \Omega^{0,1}(X,\Hom(L_j,L_i))$;

III. the Higgs field is of the form \begin{eqnarray}
\theta=\begin{pmatrix}a_{11}&a_{12}&a_{13}&\cdots&a_{1n}\\\gamma_1&a_{22}&a_{23}&\cdots&a_{2n}\\&\gamma_2&a_{33}&\cdots&a_{3n}\\&&\ddots&\ddots&\vdots\\&&&\gamma_{n-1}&a_{nn}\end{pmatrix}\label{Phi}
\end{eqnarray} where $a_{ij}\in \Omega^{1,0}(X,\Hom(L_j,L_i))$ and $\gamma_k:L_k\rightarrow L_{k+1}\otimes K$ is holomorphic.

We then consider the subbundle $F=F_k$ for $k=1,\cdots, n-1$. Then the associated factor $B$ is 
\[B=\begin{pmatrix}0&0&\cdots&0&\gamma_k\\0&0&\cdots&0&0\\\vdots&\vdots&\cdots&\vdots&\vdots\\0&0&\cdots&0&0\end{pmatrix}:F_k\rightarrow (L_{k+1}\oplus \cdots\oplus L_n)\otimes K\]
then
$$\sqrt{-1}\Lambda_{g_X} \tr(B^{*_H}\wedge B)=-|\gamma_k|^2(h|_{L_k})^{-1}h|_{L_{k+1}}/g_X=-|\gamma_k|^2\frac{\det(h|_{F_{k-1}})\det(h|_{F_{k+1}})}{\det(h|_{F_k})^2}/g_X.$$

Therefore the Hitchin equation for $(E,\theta, h)$ and $F=F_k (k=1,\cdots,n-1)$ becomes
\begin{eqnarray}
-\sqrt{-1}\Lambda_{g_X}\tr(F(h|_{F_k})) \geq -|\gamma_k|^2\frac{\det(h|_{F_{k-1}})\det(h|_{F_{k+1}})}{\det(h|_{F_k})^2}/g_X, \quad k=1,\cdots,n-1.
\end{eqnarray}
Note that the Hitchin equation for $(E_0,\theta_0, h_1)$ and $F=F_k$ 
\begin{eqnarray}
-\sqrt{-1}\Lambda_{g_X}\tr(F(h_1|_{F_k}))= -|\gamma_k|^2\frac{\det(h_1|_{F_{k-1}})\det(h_1|_{F_{k+1}})}{\det(h_1|_{F_k})^2}/g_X, \quad k=1,\cdots,n-1.
\end{eqnarray}

Set \(\displaystyle v_k=\log \frac{\det(h|_{F_k})}{\det(h_1|_{F_k})}\) for $1\leq k\leq n$ and $v_0=0$. The Laplacian with respect to $g_X$ is $2\sqrt{-1}\Lambda_{g_X}\partial\bar\partial$, denoted by $\triangle_{g_X}$. 
 We obtain
\begin{equation}\label{prekeyequation}
\frac{1}{2}\triangle_{g_X} v_k+(e^{v_{k-1}+v_{k+1}-2v_k}-1)\cdot|\gamma_k|^2\frac{\det(h_1|_{F_{k-1}})\det(h_1|_{F_{k+1}})}{\det(h_1|_{F_k})^2}/g_X\geq 0, \quad k=1,\cdots,n-1.
\end{equation}

Let $$c_k=|\gamma_k|^2\frac{\det(h_1|_{F_{k-1}})\det(h_1|_{F_{k+1}})}{\det(h_1|_{F_k})^2}/g_X\int_{0}^{1}e^{(1-t)(v_{k-1}+v_{k+1}-2v_k)}dt,\quad k=1,\cdots,n-1.$$
Then $v_k$'s satisfy
\begin{eqnarray}\label{OurSystem}
\frac{1}{2}\triangle_{g_X} v_k+c_k(v_{k-1}-2v_k+v_{k+1})&\geq&0,\quad k=1,\cdots, n-1.
\end{eqnarray}

By the assumption on the boundary $\partial X$, $v_k=0, k=1,\cdots, n-1.$
It is easy to check that the above system of equations satisfies the assumptions in Lemma \ref{MaximumPrincipleGeneral}. Moreover, $(1,1,\cdots,1)$ is indeed a supersolution of the system (\ref{OurSystem}). Then one can apply Lemma \ref{MaximumPrincipleGeneral} and obtain $v_k\leq 0, k=1,\cdots,n-1$. 
\hfill\qed

\begin{lem}(\cite[Theorem 1]{Sirakov})\label{MaximumPrincipleGeneral}
Let $(X,g)$ be a Riemannian manifold with boundary. For each $1\leq i\leq n$, let $u_i$ be a $C^2$ real-valued function on $X$ satisfying
\begin{eqnarray*}
&&\triangle_{g} u_i+\sum_{j=1}^{n}c_{ij}u_j\geq 0, \quad 1\leq i\leq n,\quad \text{in $X$}, 
\end{eqnarray*}
where $c_{ij}$ are continuous functions on $X$, $1\leq i,j\leq n$, satisfying\\
$(a)$ cooperative: $c_{ij}\geq 0,~ i\neq j$,\\
$(b)$ fully coupled: the index set $\{1,\cdots,n\}$ cannot be split up in two disjoint nonempty sets $\alpha,\beta$ such that $c_{ij}\equiv 0$ for $i\in\alpha,j\in \beta.$

Suppose that there exists a supersolution $(\psi_1,\psi_2,\cdots,\psi_n)$ satisfying $\psi_i\geq 1$ of the above system, i.e., 
\begin{eqnarray*}
&&\triangle_{g} \psi_i+\sum_{j=1}^{n}c_{ij}\psi_j\leq 0, \quad 1\leq i\leq n.
\end{eqnarray*}
Then $$\sup_{X}u_i\leq \sup_{\partial X}u_i,\quad  1\leq i\leq n.$$
\end{lem}

\subsection{Domination property and the existence of harmonic metrics}

\label{subsection;23.3.31.1}
We assume that $X$ is non-compact.

Let $(E,\theta)$, $(E_0,\theta_0)$
be as in \S\ref{Setup}. Moreover, we assume the following.
\begin{condition}
\label{condition;23.5.4.1}
There exists a harmonic metric $h_0$ of $(E_0,\theta_0)$
such that the decomposition $E_0=\oplus_{i=1}^n\Gr^F_i(E)$ is orthogonal
with respect to $h_0$.
Note that $X$ has to be hyperbolic,
see \cite[Lemma 3.13]{LiMochizuki}.
\end{condition}

Let $\Harm^{dom}(E,\theta: h_0)$ denote the set of harmonic metrics $h$ of $(E,\theta)$ such that (i) $h$ weakly dominates $h_0$, (ii) $\det(h)=\det(h_0)$. 
We shall prove the following theorem
in \S\ref{subsection;23.3.31.11}
after the preliminaries 
in \S\ref{subsection;23.3.31.10}--\S\ref{subsection;23.3.31.12}.
\begin{thm}\label{MainTheorem}
\begin{itemize}
    \item $\Harm^{dom}(E,\theta:h_0)$ is not empty.
    \item $\Harm^{dom}(E,\theta:h_0)$ is compact in the following sense: any sequence $h_i$ in $\Harm^{dom}(E,\theta:h_0)$ contains a subsequence $h_i'$ such that the sequence $h_i'$ and their derivatives are convergent on any relatively compact open subset $K$ of $X$.
\end{itemize}
\end{thm}

Let $ (E, \theta,C), (E_0,\theta_0, C)$ be defined in \S\ref{SymmetricForms}.
In addition to Condition \ref{condition;23.5.4.1},
we assume the following.
\begin{condition}
$h_0$ is compatible with $C$.
\end{condition}

Let $\Harm^{dom}(E,\theta, C: h_0)$ denote the set of harmonic metrics $h$ of $(E,\theta)$ such that (i) $h$ weakly dominates $h_0$, (ii) $h$ is compatible with $C$. 
We shall also prove the following theorem
in \S\ref{subsection;23.3.31.11}.
\begin{thm}\label{MainTheorem1}
$\Harm^{dom}(E,\theta,C:h_0)$ is non-empty and compact.
\end{thm}

\subsubsection{Preliminary from linear algebra}
\label{subsection;23.3.31.10}

Let $P$ be an upper triangular $n\times n$ matrix
with non-vanishing diagonal terms.

Let $A$ be an $n\times n$ matrix with
$$A_{j,k}=0(j>k+1),\quad A_{k+1,k}\neq 0.$$
Set
$$|A|:=\max_{j,k}|A_{j,k}|, \quad \widetilde{|A|}=\max_{1\leq k\leq n-1}|(A_{k+1,k})^{-1}|.$$  In this section, our goal is to show the following. 
\begin{prop}\label{AlgebraProposition}
Suppose $|P^{-1}AP|\leq c, |P_{1,1}|\geq d, \det(P)\leq e$,
then there exists a constant $C=C(|A|,\widetilde{|A|},c,d,e)$ such that 
$$|(P^{-1})_{i,j}|+|P_{i,j}|\leq C.$$
\end{prop}
The proof of Proposition \ref{AlgebraProposition} follows from Proposition \ref{Prop1} and \ref{Prop2}.\\

First we investigate the properties of  $P^{-1}$ in terms of $P$.
\begin{lem}\label{InverseMatrix}
\begin{itemize}
    \item $(P^{-1})_{i,j}=0$ for $i>j$.
    \item $(P^{-1})_{j,j}=(P_{j,j})^{-1}$ for $j=1,\cdots, n.$
    \item For $1\leq i<j\leq n$ and $m\in \mathbb Z_{\geq 1}$, let $\mathcal S_m(i,j)$ denote the set of $\mathbf{i}=(i_0,i_1,\cdots, i_m)\in \mathbb Z_{\geq 1}^m$ such that $i_0=i<i_1<\cdots<i_m=j$. Then, 
    $$(P^{-1})_{i,j}=\sum_{m\geq 1}\sum_{\mathbf{i}\in \mathcal S_m(i,j)}(-1)^m\prod_{p=0}^m(P_{i_p,i_p})^{-1}\prod_{p=0}^{m-1}P_{i_p,i_{p+1}}.$$
\end{itemize}
\end{lem}
\pf
Let $Q$ be the diagonal matrix
such that $Q_{i,i}=P_{i,i}$.
We set $R=P-Q$,
which is strictly upper triangular matrix.
Let $I$ denote the identity matrix.
Because $P=Q(I+Q^{-1}R)$,
we obtain
\[
 P^{-1}=
 Q^{-1}+\sum_{m\geq 1}(-1)^m(Q^{-1}R)^mQ^{-1}.
\]
Then, the claims of the lemma are obvious.
\hfill\qed

\begin{lem}\label{MutualControl}
Assume $|(P_{i,i})^{-1}|\leq B_1$.
Suppose $|P_{i,i+t}|\leq B_2$ for all $0\leq t\leq t_0$, then $$|(P^{-1})_{i,i+t}|\leq c\sum_{m=0}^{t}B_1^{m+1}B_2^{m}$$ is bounded by a constant $c=c(n)$ for all $0\leq t\leq t_0$.
\end{lem}
\pf
Note that each term of the formula of $P^{-1}_{i,i+t}$ involves terms of products of $(P_{i_p,i_p})^{-1}$ and $P_{i_p,i_{p+1}}$ for $i\leq i_p<i_{p+1}\leq j.$
By assumption, all such terms are bounded by $B$ since $i_{p+1}\leq i_p+t_0$.
\hfill\qed

\begin{prop}\label{Prop1}
Assume $B_1^{-1}\leq |P_{i,i}|\leq B_1$ and suppose $|P^{-1}AP|\leq B_2.$ Then we have 
$$|(P^{-1})_{i,j}|+|P_{i,j}|\leq C, \quad 1\leq i,j\leq n$$ for some constant $C=C(n,B_1,B_2,|A|,\widetilde{|A|}).$
\end{prop}
\pf
It is enough to estimate $P_{i,j}$ with $i\leq j.$
We prove by induction on $j-i$. First of all, $P_{i,i}$ satisfies the estimates by assumption. Assume that \begin{equation}\label{assumption}
|P_{i,i+t}|\leq C(t_0),\quad 1\leq i\leq n, 0\leq t\leq t_0.
\end{equation}
We are going to show $|P_{i,i+t_0+1}|\leq C.$

By Assumption (\ref{assumption})
and Lemma \ref{MutualControl}, $|(P^{-1})_{i,i+t}|\leq C(t_0)$ for all $i, 0\leq t\leq t_0$.
We set
\[
\nbigt(i,t_0)=\bigl\{
(\ell,k)\,\big|\,
i-1\leq \ell-1\leq k\leq i+t_0
\bigr\},
\quad
\nbigt'(i,t_0)=\nbigt(i,t_0)\setminus\{(i,i-1),(i+t_0+1,i+t_0)\}.
\]
For any $(k,\ell)\in\nbigt'(i,t_0)$,
we have
$\ell-i\leq t_0$ and $i+t_0-k\leq t_0$.
We obtain
\begin{eqnarray*}
&&(P^{-1}AP)_{i,i+t_0}\\
&=&\sum_{(l,k)\in\nbigt(i,t_0)}(P^{-1})_{i,l}A_{l,k}P_{k,i+t_0}\\
&=&\sum_{(l,k)\in\nbigt'(i,t_0)}(P^{-1})_{i,l}A_{l,k}P_{k,i+t_0}
+(P^{-1})_{i,i+t_0+1}A_{i+t_0+1,i+t_0}P_{i+t_0,i+t_0}
+(P^{-1})_{i,i}A_{i,i-1}P_{i-1,i+t_0} \\
&=&
\sum_{(l,k)\in\nbigt'(i,t_0)}(P^{-1})_{i,l}A_{l,k}P_{k,i+t_0}
+\sum_{m\geq 2}
\sum_{\mathbf{i}\in\nbigs_m(i,i+t_0+1)}
(-1)^m
 A_{i+t_0+1,i+t_0}P_{i+t_0,i+t_0}
 \prod_{p=0}^{m}(P_{i_0,i_0})^{-1}
 \prod_{p=0}^{m-1}P_{i_p,i_{p+1}}
\\
 & &
  -P_{i,i+t_0+1}P_{i,i}^{-1}P_{i+t_0+1,i+t_0+1}^{-1}
  A_{i+t_0+1,i+t_0}P_{i+t_0,i+t_0}
+(P^{-1})_{i,i}A_{i,i-1}P_{i-1,i+t_0}.
\end{eqnarray*}
Here, we formally put
$A_{1,0}=A_{n+1,n}=0$,
$(P^{-1})_{i,n+1}=P_{i,n+1}=P_{0,1+t_0}=0$
and $P_{n+1,n+1}=1$.
The first and second terms
of the right hand side in the formula of $(P^{-1}AP)_{i,i+t_0}$
only involve $P_{\mu\nu}$ where $\nu\leq \mu+t_0$ and $A_{lk}$.
By the formula in the case $i=1$,
we obtain an estimate for
$P_{1,t_0+2}$.
Inductively,
we obtain an estimate for $P_{i,i+t_0+1}$ $(i=2,\ldots,n-t_0-1)$
by using the formula in the case $i$.
\hfill\qed

\begin{prop}\label{Prop2}
 Suppose $|(P^{-1}AP)_{k+1,k}|\leq c$,
$P_{1,1}\geq d$ and $\det(P)\leq e$.
Then 
\begin{equation}
\label{eq;23.4.2.1}
 d\bigl(
 c\widetilde{|A|}
 \bigr)^{1-i}
 \leq
 |P_{i,i}|
 \leq
 (ed^{-i+1})^{\frac{1}{n+1-i}}
 (c\widetilde{|A|})^{-\frac{(i-1)(i-2)}{2(n+1-i)}}
 (c\widetilde{|A|})^{-\frac{1}{2}(n-i)}
\end{equation}
\end{prop}
\pf
We set $\ctilde=c\widetilde{|A|}$ to simplify the notation.
Recall that $(P^{-1}AP)_{k+1,k}=(P^{-1})_{k+1,k+1}A_{k+1,k}P_{k,k}.$
Thus $$|(P^{-1})_{k+1,k+1}P_{k,k}|\leq \ctilde.$$
By Lemma \ref{InverseMatrix}, $(P^{-1})_{k+1,k+1}=P_{k+1,k+1}^{-1}$.
Thus
\[
 |P_{k,k}|\leq \ctilde|P_{k+1,k+1}|.
\]
So
\[
 |P_{j,j}|\geq
 \ctilde^{j-1}|P_{1,1}|
 \geq
 \ctilde^{j-1}d.
\]

Because
$|P_{j,j}|\geq \ctilde^{j-i}|P_{i,i}|$ for $i\leq j$,
we obtain
\[
 \prod_{j=1}^{i-1}(\ctilde^{j-1}d)
 \cdot
 \prod_{j=i}^n(\ctilde^{j-i}|P_{i,i}|)
 \leq
 \prod_{j=1}^n|P_{j,j}|
 =e.
\]
It implies
\[
 |P_{i,i}|^{n+1-i}
 \leq
 ed^{-i+1}\ctilde^{-\frac{1}{2}(i-1)(i-2)}
 \cdot
 \ctilde^{-\frac{1}{2}(n-i)(n+1-i)}.
\]
Thus, we obtain the right inequality
in (\ref{eq;23.4.2.1}).
\hfill\qed

\subsubsection{Notation}

Let $V$ be a complex vector space equipped with a base $\mathbf{e}=(e_1,\cdots, e_n)$. Let $h$ be any Hermitian metric of $V$. By applying the Gram-Schmidt process to the base $\mathbf{e}$, we obtain a base $\mathbf{v}(h)=(v_1(h), \cdots, v_n(h)).$ Let $P(h)=(P(h)_{j,k})$ be the matrix determined by $\mathbf{v}=\mathbf{e}P(h)$. Then $P(h)_{j,k}=0(j>k)$, i.e., 
$$v_k(h)=\sum_{j\leq k}P(h)_{j,k}e_j.$$
Let $P^{-1}(h)=(P^{-1}(h)_{j,k})$ be the inverse matrix of $P(h)$.
In terms of the frame $\vece$,
the metric $h$ is represented by the matrix
$h(\vece)=(P(h)^{-1})^t\cdot \overline{P(h)^{-1}}$.

We use a similar notation
for a vector bundle equipped with a frame and a Hermitian metric.

\subsubsection{Local estimate in the nowhere vanishing case}

We set $U(R)=\{z\in\mathbb C||z|<R\}$
and $\overline U(R)=\{z\in\mathbb C||z|\leq R\}$ for any $R>0$.

Let $R_1<R_2$. Let $E=\oplus_{i=1}^n\mathcal O_{U(R_2)}e_i$. Let $f$ be an endomorphism of $E$. Let $A$ be the matrix determined by $f(e_j)=\sum_{i=1}^n A_{ij}e_i.$ That is $A$ is the matrix representation of $f$ in terms of $\vece.$ Note that $|f|_h=|P(h)^{-1}A P(h)|.$ We assume the following. 
\begin{condition}
\begin{itemize}
    \item $A_{ij}=0(i>j+1)$.
    \item $A_{j+1,j}(j=1,\cdots, n-1)$ are nowhere vanishing on $\overline U(R_1).$
    \item $|\tr(f^l)|(l=1,\cdots, n)$ are bounded on $U(R_2).$
\end{itemize}
\end{condition}
We set 
\[B_1(f)=\max_{1\leq l\leq n}\sup_{U(R_2)}|\tr(f^l)|\]
\[B_2(f)=\min_{1\leq j\leq n-1}\min_{\overline U(R_1)}|A_{j+1,j}|>0.\]
We obtain the Higgs field $\theta=fdz$ of $E$. We recall the following lemma. 
\begin{lem}(\cite[Proposition 3.12]{LiMochizuki})\label{BoundedDifferential}
There exists $C_1>0$ depending only on $R_1,R_2, n$ and $B_1(f)$ such that $|f|_h\leq C_1$ on $U(R_1)$ for any harmonic metric $h$ of $(E,\theta)$ on $U(R_2).$
\end{lem}

Let $f_0$ be the endomorphism of $E$ determined by $f_0(e_j)=A_{j+1,j}e_{j+1}$ for $j=1,\cdots, n-1$ and $f_0(e_n)=0.$ We obtain the Higgs field $\theta_0=f_0dz$ of $E$. Assume that there exists a harmonic metric $h_0$ of $(E,\theta_0)$ such that the decomposition $E=\oplus_{i=1}^n \mathcal O_{U(R_2)}e_i$ is orthogonal. 

Let $\Harm^{dom}(E,\theta:h_0)$ denote the set of harmonic metrics $h$ of $(E,\theta)$ such that (i) $h$ weakly dominates $h_0$, (ii) $\det(h)=\det(h_0)$. For two Hermitian metrics $h_1,h_2$ on $E$, let $s(h_1,h_2)$ be the automorphism of $E$ such that $$h_2(u,v)=h_1(s(h_1,h_2)u, v),$$ for any two sections $u,v$
of $E$. In terms of the frame $\vece$, $s(h_1,h_2)$ is represented by a matrix $J(h_1,h_2)$.
Then $J(h_1,h_2)$ satisfies $h_2(\vece)=J(h_1,h_2)^t\cdot h_1(\vece).$ 
So $$J(h_1,h_2)=P(h_1)\cdot \overline{P(h_1)^t}\cdot\overline{(P(h_2)^{-1})^t}\cdot P(h_2)^{-1}.$$

We obtain the following proposition.
\begin{prop}\label{InitialControl}
There exists $C_2>0$ depending only on $n, R_i(i=1,2)$ and $B_k(f)(k=1,2)$ such that 
\[|s(h_0,h)|_{h_0}+|s(h_0,h)^{-1}|_{h_0}\leq C_2\]
for any $h\in \Harm^{dom}(E,\theta:h_0).$
\end{prop}
\pf
From Lemma \ref{BoundedDifferential}, we have $$|f|_h=|P(h)^{-1}A P(h)|\leq C_1.$$ 

Since $h$ weakly dominates $h_0,$ we obtain that $|P(h)_{1,1})\geq |P(h_0)_{1,1}|\geq C_1'$ for some positive constant $C_1'$. And $\det(P(h))=\det(P(h_0))\leq C_2',$ for some positive constant $C_2'$.

From Proposition \ref{AlgebraProposition}, we have $|P(h)|+|P(h)^{-1}|\leq C_3'$, for some positive constant $C_3'$.

The rest follows from the matrix expression $J(h_0,h)$ of $s(h_0,h)$ is
$$J(h_0,h)=P(h_0)\cdot\overline{P(h_0)^t}\cdot\overline{(P(h)^{-1})^t}\cdot P(h)^{-1},$$
$$|s(h_0,h)|_{h_0}=|P(h_0)^{-1}J(h_0,h)P(h_0)|$$ and $P(h_0),P(h_0)^{-1}$ are bounded.
\hfill\qed

\subsubsection{Local estimate in the general case}
\label{subsection;23.3.31.12}

Let $X$, $(E,\theta)$, $(E_0,\theta_0)$ and $h_0$ be
as in \S\ref{subsection;23.3.31.1}.
We fix an isomorphism
$E\simeq E_0$ as in \S\ref{subsection;23.3.31.2},
and we regard $h_0$ as a Hermitian metric of $E$.

Let $K_1\subset X$ be a relatively compact open subset.
Let $K_2$ be a relatively compact open neighbourhood of $\overline{K}_1$ in $X$.
\begin{prop}\label{BoundedProposition}
There exists $C_3>0$ such that the following holds on $K_1$ for any $h\in \Harm^{dom}((E,\theta:h_0)|_{K_2})$:
\begin{equation}
    \label{DesiredControl}
|s(h_0,h)|_{h_0}+|s(h_0,h)^{-1}|_{h_0}\leq C_3.
\end{equation}
\end{prop}
\pf
By making $K_1$ larger if necessary, we may assume that $A_{j+1,j}(j=1,\cdots, r-1)$ are nowhere vanishing on a neighbourhood $N$ of $\partial K_1$. Let $N'$ be a relatively compact neighbourhood of $\partial K_1$ in $N$. By using Proposition \ref{InitialControl}, we can prove that there exists $C_4>0$ such that the following holds on $N'$ for any $h\in \Harm^{dom}((E,\theta:h_0)|_{N})$:
\begin{equation}
\label{Control1}
|s(h_0,h)|_{h_0}+|s(h_0,h)^{-1}|_{h_0}\leq C_4.
\end{equation}
Let $h_1$ be a harmonic metric of $(E,\theta)|_{\overline K_2}$ such that $\det(h_1)=\det(h_0)$. There exists $C_5>0$ such that the following holds on $\overline K_2$:
\begin{equation}
\label{Control2}
|s(h_1,h_0)|_{h_1}+|s(h_1,h_0)^{-1}|_{h_1}\leq C_5.
\end{equation}
By Equation (\ref{Control1}) and (\ref{Control2}), there exists $C_6>0$ such that the following holds for any $h\in \Harm^{dom}((E,\theta:h_0)|_{K_2})$ on $N':$ 
\[|s(h_1,h)|_{h_1}+|s(h_1,h)^{-1}|_{h_1}\leq C_6.\]
Because $\log \tr(s(h_1,h))$ are subharmonic, the following holds on $K_1:$
\[|s(h_1,h)|_{h_1}+|s(h_1,h)^{-1}|_{h_1}\leq C_6.\]
Therefore, together with Equation (\ref{Control2}), we obtain Equation (\ref{DesiredControl}). 
\hfill\qed

\subsubsection{Proof of Theorem \ref{MainTheorem} and Theorem \ref{MainTheorem1}}
\label{subsection;23.3.31.11}
Let $X_i(i=1,2,\cdots)$ be a smooth exhaustion family of $X$. Let $h^{(i)}$ be the harmonic metrics of $(E,\theta)|_{X_i}$ such that $h^{(i)}|_{\partial X_i}=h_0|_{\partial X_i}$. 

\begin{thm}\label{Convergence}
$h^{(i)}$ contains a convergent subsequence.
\end{thm}
\pf 
By Proposition \ref{Control}, $h^{(i)}$ weakly dominates
$h_0$. By Proposition \ref{prop;20.5.29.1} and Proposition \ref{BoundedProposition}, $h^{(i)}$ contains a convergent subsequence.
\hfill\qed

\vspace{.1in}
Hence, we obtain the first claim of Theorem \ref{MainTheorem}. We also obtain the second claim of Theorem \ref{MainTheorem} from Proposition \ref{prop;20.5.29.1} and the argument in the proof of Proposition \ref{prop;20.5.29.1}.

\vspace{.1in}
Suppose moreover, $E, \theta, C, E_0,\theta_0, h_0$ are in the setting of Theorem \ref{MainTheorem1}. By the uniqueness of solutions to the Dirichlet problem and $h_0$ is compatible with $C$, $h^{(i)}$ is also compatible with $C$. So the limit metric is again compatible with $C$. So we obtain the first claim of Theorem \ref{MainTheorem1}. The second claim of Theorem \ref{MainTheorem1} follows from the second claim of Theorem \ref{MainTheorem} and that compatibility with $C$ is preserved under limit.


\section{Uniqueness in a bounded case}
\label{UniquenessSection}

\subsection{Statement}

Let $X$ be a Riemann surface.
Let $g_X$ be a complete K\"ahler metric 
whose Gauss curvature is bounded below.

We fix a line bundle $K_X^{1/2}$ and an isomorphism
$K_X^{1/2}\otimes K_X^{1/2}\simeq K_X$.
We set
\[
 \hyperk_{X,n}=\bigoplus_{i=1}^n K_X^{(n+1-2i)/2}.
\]
We set
$F_j\hyperk_{X,n}=\bigoplus_{i\leq j}K_X^{(n+1-2i)/2}$.
We obtain the Hermitian metric
$h_X=\bigoplus g_{X}^{-(n+1-2i)/2}$ of $\hyperk_{X,n}$.
Let $C$ be a holomorphic non-degenerate symmetric pairing
of $\hyperk_{X,n}$ which is compatible with $h_X$.

Let $\theta$ be a Higgs field of $\hyperk_{X,n}$.
We assume the following.
\begin{itemize}
 \item $\theta (F_j\hyperk_{X,n})\subset F_{j+1}\hyperk_{X,n}
       \otimes K_X$.
       Moreover, the induced morphisms
       $\phi_j\colon
       \Gr^F_j\hyperk_{X,n}\to
       \Gr^F_{j+1}\hyperk_{X,n}\otimes K_X$
       are the identity morphisms
       under
       the natural isomorphisms
       $\Gr^F_j\hyperk_{X,n}=
       K_X^{(n+1-2j)/2}
       =\Gr^F_{j+1}\hyperk_{X,n}\otimes K_X$.
 \item $\theta$ is bounded
       with respect to $h_X$ and $g_X$.
 \item $\theta$ is self-adjoint
       with respect to $C$.
\end{itemize}

We shall prove the uniqueness of harmonic metrics
which are compatible with $C$
and mutually bounded with $h_X$.

\begin{thm}
\label{thm;23.4.2.11}
Let $h_1$ and $h_2$ be harmonic metrics of 
$(\hyperk_{X,n},\theta)$.
Suppose that
both $h_i$ are compatible with $C$,
and that both $h_i$ are mutually bounded with $h_X$.
Then, $h_1=h_2$ holds.
\end{thm}

\subsubsection{A characterization of the mutual boundedness with $h_X$}

We also have the following characterization
for a harmonic metric to be mutually bounded with $h_X$.
\begin{prop}
\label{prop;23.4.2.10}
Let $h$ be a harmonic metric of $(\hyperk_{X,n},\theta)$
such that $\det(h)=1$.
Then, $h$ is mutually bounded with $h_X$
if and only there exists $b>0$ such that
$h_{|F_1\hyperk_{X,n}}\leq b h_{X|F_1\hyperk_{X,n}}$.
\end{prop}
\pf
The ``only if'' part of Proposition \ref{prop;23.4.2.10} is clear.
Let us prove the ``if'' part.
Let $h$ be a harmonic metric of $(\hyperk_{X,n},\theta)$
such that $\det(h)=1$.
Because the spectral curve of
the Higgs bundle $(\hyperk_{X,n},\theta)$ is bounded
with respect to $g_X$,
we obtain the following lemma from 
\cite[Proposition 3.12]{LiMochizuki}.
\begin{lem}
$|\theta|_{h,g_X}$ is bounded on $X$.
\hfill\qed 
\end{lem}

Let $x$ be any point of $X$.
Let $\tau$ be a base of $K_{X|x}$
such that $|\tau|_{g_{X|x}}=1$.
By setting
$e_i=\tau^{(n+1-2i)/2}$ $(i=1,\ldots,n)$,
we obtain an orthonormal frame
$\vece=(e_1,\ldots,e_n)$ of $\hyperk_{X,n|x}$
with respect to $h_{X|x}$.
Let $A$ be the matrix determined by
$\theta_{|x}(\vece)=\vece\cdot A\,\tau$.
Because $\theta$ is bounded with respect to $h_X$ and $g_X$,
there exists $B_1>0$ which is independent of $x$
such that $|A|\leq B_1$.
Moreover,
$A_{k+1,k}=1$ for $k=1,\ldots,n-1$.

By applying the Gram-Schmidt process to
the frame $\vece$ and the metric $h_{|x}$,
we obtain the base $\vecv$ of $\hyperk_{X,n|x}$
which is orthonormal with respect to $h_{|x}$.
Let $P$ be the matrix determined by
$\vecv=\vece\cdot P$.
Because $\theta$ is bounded with respect to $h$ and $g_X$,
there exists $B_2>0$, which is independent of $x$,
such that $|P^{-1}AP|\leq B_2$.
Because $h_{|F_1\hyperk_{X,n}}\leq b h_{X|F_1\hyperk_{X,n}}$,
we obtain $P_{1,1}\geq b^{-1}$.
Because $\det(h)=1$,
we have $\det(P)=1$.
By Proposition \ref{AlgebraProposition},
there exists $B_3>0$ which is independent of $x$
such that
$|P|+|P^{-1}|\leq B_3$.
Therefore, there exists $B_4$
such that the following holds on $X$:
\[
 \bigl|s(h,h_X)\bigr|_{h_X}
 +\bigl|s(h,h_X)^{-1}\bigr|_{h_X}
 \leq B_4
\]
Thus, we obtain Proposition \ref{prop;23.4.2.10}.
\hfill\qed

\subsection{Preliminary from Linear algebra}

\subsubsection{Cyclic vectors}

Let $V$ be an $n$-dimensional complex vector space
equipped with an endomorphism $f$.
A vector $v\in V$ is called an $f$-cyclic vector
if $v,f(v),\ldots,f^{n-1}(v)$ generate $V$.
The following proposition is well known.
(For example, see \cite[\S6,\S7]{Linear-Algebra-Textbook}.)
\begin{prop}
There exists an $f$-cyclic vector
if and only if
the characteristic polynomial of $f$
equals the minimal polynomial of $f$.
\hfill\qed
\end{prop}

\begin{cor}
\label{cor;22.9.12.3}
Suppose that there exists an $f$-cyclic vector.
For any eigenvalue $\alpha$ of $f$,
the space of eigen vectors associated with $\alpha$
is one dimensional. 
\hfill\qed
\end{cor}

Let $h$ be a Hermitian metric of $V$.
For any $v\in V$,
we set $\omega(f,v)=v\wedge f(v)\wedge\cdots\wedge f^{n-1}(v)$.
Then, $v$ is an $f$-cyclic vector of $V$
if and only if $\omega(f,v)\neq 0$.
We always have
$|\omega(f,v)|_h\leq |f|_h^{n(n-1)/2}|v|^n_h$.

\begin{lem}
\label{lem;23.2.22.20}
Let $A>0$ and $\rho>0$.
There exists $\epsilon_0(n,A,\rho)>0$
depending only on $n$, $A$ and $\rho$,
such that the following holds.
\begin{itemize}
 \item Suppose that
       $|f|_h\leq A$
       and that 
       $\rho|v|_h^n\leq |\omega(f,v)|_h$
       for a non-zero element $v\in V$.
       Let $f_1$ be an endomorphism of $V$
       such that
       $|f-f_1|_h\leq \epsilon_0(n,A,\rho)$.
       Then,
       $\frac{1}{2}\rho|v|_h^n<|\omega(f_1,v)|_h$.
       In particular,
       $f_1$ also has a cyclic vector.
\end{itemize}  
\end{lem}
\pf
If $|f-f_1|_h<1$,
we obtain
\[
\bigl|
f_1^j(v)-f^j(v)
\bigr|_h
\leq
\sum_{k=1}^jC(j,k)|f|^{j-k}_h|f-f_1|^k_h|v|_h
\leq |f-f_1|_h(1+|f|_h)^j|v|_h.
\]
Here, $C(j,k)$ denote the binomial coefficients.
We obtain
\begin{multline}
\Bigl|
v\wedge f_1(v)\wedge\cdots \wedge f_1^{j-1}(v)
\wedge(f_1^{j}(v)-f^j(v))
\wedge f^{j+1}(v)\wedge\cdots f^{n-1}(v)
\Bigr|_h
\leq
|v|_h^n\cdot |f_1|_h^{j(j-1)/2}\cdot |f^j-f_1^j|_h
\cdot |f|_h^{n(n-1)/2-j(j+1)/2}
 \\
 \leq
|v|_h^n\cdot
 |f-f_1|_h
 \cdot (1+|f|_h)^{n(n-1)/2}.
\end{multline}
We obtain
\[
 \bigl|\omega(f,v)-\omega(f_1,v)\bigr|_h
 \leq
 n(1+|f|_h)^{n(n-1)/2}|f-f_1|_h\cdot |v|_h^{n}.
\]
Then, the claim of the lemma is clear.
\hfill\qed

\subsubsection{Real structure and self-adjoint endomorphisms}
\label{subsection;22.9.12.2}

Let $C$ be a non-degenerate symmetric pairing of
a finite dimensional complex vector space $V$.
Let $f$ be an endomorphism of $V$
such that $f$ is self-adjoint with respect to $C$,
i.e., $C(fu,v)=C(u,fv)$ for any $u,v\in V$.
There exists the generalized eigen decomposition
$V=\bigoplus_{\alpha\in\cnum} V_{\alpha}$,
where $V_{\alpha}$ denote the space of
generalized eigen vectors of $f$ associated with $\alpha$.
The following lemma is well known.

\begin{lem}
\label{lem;23.2.22.1}
If $\alpha\neq\beta$,
then $V_{\alpha}$ and $V_{\beta}$ are orthogonal
with respect to $C$. 
\end{lem}
\pf
We explain a proof just for the convenience of readers.
For $j\in\seisuu_{\geq 1}$,
we set
\[
 \nbigf_jV_{\alpha}
 =\bigl\{
 u_{\alpha}\in V_{\alpha}\,\big|\,
 (f-\alpha\id_V)^ju_{\alpha}=0
 \bigr\}.
\]
Let us prove that
$\nbigf_iV_{\alpha}$ and
$\nbigf_jV_{\beta}$ $(i+j= \ell)$
are orthogonal 
by an induction on $\ell$.
Let us consider the case $\ell=2$.
For $u_{\alpha}\in\nbigf_1V_{\alpha}$
and $v_{\beta}\in\nbigf_1V_{\beta}$,
we obtain
$\alpha C(u_{\alpha},v_{\beta})
=C(f(u_{\alpha}),v_{\beta}) 
=C(u_{\alpha},f(v_{\beta}))
=\beta C(u_{\alpha},v_{\beta})$.
It implies $C(u_{\alpha},v_{\beta})=0$.
Suppose that we have already proved the claim
in the case $i+j=\ell$,
and let us consider the case $i+j=\ell+1$.
For $u_{\alpha}\in\nbigf_iV_{\alpha}$
and $v_{\beta}\in\nbigf_jV_{\beta}$,
we have
$f(u_{\alpha})-\alpha u_{\alpha}\in \nbigf_{i-1}V_{\alpha}$
and
$f(u_{\beta})-\beta u_{\beta}\in \nbigf_{j-1}V_{\beta}$.
By the assumption of the induction,
we obtain
$C(f(u_{\alpha})-\alpha u_{\alpha},v_{\beta})=0$
and 
$C(u_{\alpha},f(v_{\beta})-\beta v_{\beta})=0$.
Because $C(f(u_{\alpha}),v_{\beta})=C(u_{\alpha},f(v_{\beta}))$,
we obtain $C(u_{\alpha},v_{\beta})=0$.
\hfill\qed

\vspace{.1in}

Let $h$ be a Hermitian metric compatible with $C$.
Let $\kappa$ be the real structure of $V$ induced by $C$ and $h$.
Let $W\subset V$ be a vector subspace
such that
(i) $f(W)\subset W$ and $f(\kappa(W))\subset\kappa(W)$,
(ii) $W\cap\kappa(W)=0$.

\begin{prop}
 We have either
(i) $W=0$,
or (ii) $f_{|W}$ and $f_{|\kappa(W)}$
have a common eigenvalue. 
\end{prop}
\pf
Suppose that there is no common eigenvalue of
$f_{|W}$ and $f_{|\kappa(W)}$.
By Lemma \ref{lem;23.2.22.1},
$W$ and $\kappa(W)$ are orthogonal with respect to $C$.
For any $u\in W$,
we obtain
$h(u,u)=C(u,\kappa(u))=0$,
and hence $W=0$.
\hfill\qed

\begin{cor}
\label{cor;22.9.12.4}
Suppose moreover that
there exists an $f$-cyclic vector.
Then, we obtain $W=0$. 
\end{cor}
\pf
If $W\neq 0$,
then $f_{|W}$ and $f_{|\kappa(W)}$
have a common eigenvalue $\alpha$.
Hence,
the dimension of the eigen space associated with $\alpha$
is larger than $2$,
which contradicts Corollary \ref{cor;22.9.12.3}.
\hfill\qed

\vspace{.1in}

Let us explain how to use Corollary \ref{cor;22.9.12.4}
in a simple case.

\begin{prop}
\label{prop;23.2.22.30}
Let $s$ be an automorphism of $V$
such that
(i) $s$ is Hermitian and positive definite with respect to $h$,
(ii) $s$ is self-adjoint with respect to $C$,
(iii) $[f,s]=0$.
Suppose that there exists an $f$-cyclic vector.
Then, we obtain $s=\id_V$.  
\end{prop}
\pf
There exists the eigen decomposition
$V=\bigoplus_{a>0} V_a$ of $s$.
Because $[s,f]=0$,
we obtain $f(V_a)\subset V_a$ for any $a$.
Recall $\kappa(V_a)=V_{a^{-1}}$
as in \cite[Lemma 2.10]{LiMochizukiGeneric}.
Hence, we obtain $V_a=0$ for any $a\neq 1$
by Corollary \ref{cor;22.9.12.4}.
\hfill\qed

\begin{rem}
Theorem {\rm\ref{thm;22.9.12.10}} below
is a quantified version of Proposition {\rm\ref{prop;23.2.22.30}}.
\hfill\qed
\end{rem}

\subsection{An estimate}

Let $V$ be a complex vector space
equipped with a base $\vece=(e_1,\ldots,e_n)$.
Let $C$ be a non-degenerate symmetric pairing of $V$.
For $\rho>0$,
let $\nbigh(V,\vece,C;\rho)$ be the space of Hermitian metrics $h$
of $V$ such that
(i) $h$ are compatible with $C$,
(ii) $|e_1\wedge\cdots \wedge e_n|_h\geq \rho|e_1|^n$.

Let $f$ be an endomorphism of $V$
which is self-adjoint with respect to $C$.
Let $\nbiga(f)=(\nbiga(f)_{i,j})$ be the matrix representing $f$
with respect to $\vece$,
i.e.,
$f(e_k)=\sum_{j=1}^n\nbiga(f)_{j,k}e_j$.
We assume that
$\nbiga(f)_{j,k}=0$ $(j>k+1)$
and $\nbiga(f)_{k+1,k}=1$,
i.e.,
\[
 f(e_k)=e_{k+1}+\sum_{j\leq k}\nbiga(f)_{j,k}e_j
 \quad (k=1,\ldots,n-1),
 \quad
 f(e_n)=\sum_{j\leq n}\nbiga(f)_{j,n}e_j.
\]

\begin{thm}
\label{thm;22.9.12.10}
Let $A>0$ and $\rho>0$.
There exist $\epsilon_1(n,A,\rho)>0$ and $C_1(n,A,\rho)>0$
depending only on $n$, $A$ and $\rho$
such that the following holds
for any $0<\epsilon<\epsilon_1(n,A,\rho)$:
 \begin{itemize}
  \item 
	Suppose $|f|_h\leq A$.
	Then, for any $h,h'\in\nbigh(V,\vece,C;\rho)$
	such that $\bigl|[s(h,h'),f]\bigr|_h\leq \epsilon$,
	we obtain
\[
	|s(h,h')-\id_V|_h\leq C_1(n,A,\rho)\epsilon.
\]
 \end{itemize} 
\end{thm}
\pf
Let $h,h'\in\nbigh(V,\vece,C;\rho)$.
We obtain the automorphism $s(h,h')$ of $V$
determined by
$h'(u,v)=h(s(h,h')u,v)$ for any $u,v\in V$,
which is self-adjoint with respect to both
$h$ and $h'$.
There exists the eigen decomposition
$V=\bigoplus_{a>0} V_a$ of $s(h,h')$.

Let $\kappa$ be the real structure induced by $C$ and $h$.
Note that $\kappa(V_a)=V_{a^{-1}}$.
We set $\nbigs(h,h'):=\{a> 1\,|\,V_a\neq 0\}$.
If $\nbigs(h,h')=\emptyset$,
we obtain $s(h,h')=\id_V$.
Let us consider the case where
$\nbigs(h,h')\neq\emptyset$.

Let $\nu_1$ be any positive number
such that
$\nu_1\leq\min\{1,\max\nbigs(h,h')-1\}$.
Let $c_1<c_2<\cdots<c_m$ denote the elements of
$\nbigs(h,h')$.
We set $c_0=1$.
Because $|\nbigs(h,h')|\leq n/2$,
there exists $1\leq m(0)\leq m$ such that the following holds.
\begin{itemize}
 \item $c_{i}-c_{i-1}\leq \frac{1}{2}n^{-1}\nu_1$
       for any $i<m(0)$.
 \item $c_{m(0)}-c_{m(0)-1}>\frac{1}{2}n^{-1}\nu_1$.
\end{itemize}
We set $\nbigs(h,h';\nu_1)_0=\{c_1,\ldots,c_{m(0)-1}\}$
and $\nbigs(h,h';\nu_1)_1=\{c_{m(0)},\ldots,c_{m}\}$.

\begin{lem}
\label{lem;23.2.22.10}
The set $\nbigs(h,h';\nu_1)_0$ is contained in $\{1<a\leq 2\}$.
The set $\nbigs(h,h';\nu_1)_1$ is non-empty.
For any $a_0\in\nbigs(h,h';\nu_1)_0\cup\{1\}$
and $a_1\in\nbigs(h,h';\nu_1)_1$,
we obtain
$|a_0^{-1}-a_1^{-1}|\geq \frac{1}{12}n^{-1}\nu_1$.
\end{lem}
\pf
Because $\nu_1\leq 1$, we obtain the first claim.
The second claim is clear.
For any $a_0\in\nbigs(h,h';\nu_1)_0\cup\{1\}$
and $a_1\in\nbigs(h,h';\nu_1)_1$,
we obtain
\[
|a_0^{-1}-a_1^{-1}|
\geq
\bigl|a_0^{-1}-(a_0+n^{-1}\nu_1/2)^{-1}\bigr|
=|a_0|^{-1}|a_0+n^{-1}\nu_1/2|^{-1}\frac{1}{2}n^{-1}\nu_1
\geq \frac{1}{2}\cdot \frac{1}{3}\cdot \frac{1}{2}n^{-1}\nu_1
=\frac{1}{12}n^{-1}\nu_1.
\]
Thus, we obtain the third claim of
Lemma \ref{lem;23.2.22.10}.
\hfill\qed

\vspace{.1in}

We set
\[
 W^{(\nu_1)}=\bigoplus_{a\in\nbigs(h,h';\nu_1)_1}V_a,
 \quad\quad
 V^{(\nu_1)}=
 V_1\oplus
 \bigoplus_{a\in\nbigs(h,h';\nu_1)_0}V_a
 \oplus
 \bigoplus_{a^{-1}\in\nbigs(h,h':\nu_1)_0}V_a. 
\]
Because $\nbigs(h,h';\nu_1)_1\neq\emptyset$,
we have $W^{(\nu_1)}\neq 0$.
We have
$W^{(\nu_1)}\cap \kappa(W^{(\nu_1)})=0$
and the decomposition
\[
V=V^{(\nu_1)}\oplus W^{(\nu_1)}\oplus \kappa(W^{(\nu_1)}).
\]
We obtain the decomposition
\[
 f=\sum_{U_1,U_2=V^{(\nu_1)},W^{(\nu_1)},\kappa(W^{(\nu_1)})}
 f_{U_1,U_2},
\]
where $f_{U_1,U_2}\in\Hom(U_2,U_1)$.
We set
\[
 \ftilde^{(\nu_1)}=f_{V^{(\nu_1)},V^{(\nu_1)}}
 +f_{W^{(\nu_1)},W^{(\nu_1)}}
 +f_{\kappa(W^{(\nu_1)}),\kappa(W^{(\nu_1)})}.
\]
\begin{lem}
\label{lem;23.2.22.11}
$\ftilde^{(\nu_1)}$ is self-adjoint with respect to $C$.
\end{lem}
\pf
To simplify the description,
we denote $W^{(\nu_1)}$ by $W$.
We set $\Wtilde=W\oplus \kappa(W)$.
The decomposition
$V^{(\nu_1)}\oplus \Wtilde$
is orthogonal with respect to $C$.
We obtain the decomposition
$f=\sum_{U_1,U_2=V^{(\nu_1)},\Wtilde} f_{U_2,U_1}$.
Because $f$ is self-adjoint with respect to $C$,
we obtain that
$f_{V^{(\nu_1)},V^{(\nu_1)}}$
and
$f_{\Wtilde,\Wtilde}$
are self-adjoint with respect to $C$.

We have the decompositions
$\Wtilde=W\oplus\kappa(W)$
and
$f_{\Wtilde,\Wtilde}
=\sum_{U_1,U_2=W,\kappa(W)}f_{U_2,U_1}$.
The restrictions of $C$
to $W$ and $\kappa(W)$ are $0$.
Then, it is easy to check
that
$f_{W,W}
+f_{\kappa(W),\kappa(W)}$
is self-adjoint with respect to $C$.
Thus, we obtain Lemma \ref{lem;23.2.22.11}.
\hfill\qed

\vspace{.1in}

\begin{lem}
\label{lem;23.2.22.12}
We have
 $|f-\ftilde^{(\nu_1)}|_{h}
 \leq
 \nu_1^{-1}(10n)^3
 \bigl|
 [f,s(h,h')]
 \bigr|_h$.
\end{lem}
\pf
We denote $s(h,h')$ by $s$
to simplify the description.
We have the decomposition
\[
 [s,f]=\sum_{U_1,U_2=V^{(\nu_1)},W,\kappa(W)}
 [s,f_{U_1,U_2}].
\]
We have
$\bigl|
[s,f_{U_1,U_2}]
\bigr|_h
\leq\bigl|[s,f]\bigr|_h$.

Let $U_1\neq U_2$.
Let $F_{U_2,U_1}:\Hom(U_1,U_2)\to\Hom(U_1,U_2)$
be defined by
\[
F_{U_2,U_1}(g)=[s,g]=s_{|U_2}\circ g-g\circ s_{|U_1}.
\]
For any eigenvalues $a_i$ $(i=1,2)$ of $s_{|U_i}$,
we have $|a_1-a_2|>(12 n)^{-1}\nu_1$.
Hence, $F_{U_2,U_1}$ is invertible,
and $|F_{U_2,U_1}^{-1}|_h\leq \nu_1^{-1}(12n)n^2$.
Thus, we obtain Lemma \ref{lem;23.2.22.12}.
\hfill\qed

\vspace{.1in}
By using a positive constant
$\epsilon_0(n,A,\rho)$ 
in Lemma \ref{lem;23.2.22.20},
we set
\[
\epsilon_1(n,A,\rho)
:=
\frac{1}{2}(10n)^{-3}
\epsilon_0(n,A,\rho),
\quad
C_1(n,A,\rho):=n\epsilon_1(n,A,\rho)^{-1}.
\]

Let $0<\epsilon<\epsilon_1(n,A,\rho)$.
Suppose $\bigl|[s(h,h'),f]\bigr|_h\leq \epsilon$.
We set
\[
\nu_2:=\frac{1}{2}\epsilon_1(n,A,\rho)^{-1}\epsilon
< \frac{1}{2}.
\]
If $\nu_2\leq\max\nbigs(h,h')-1$,
we obtain
\[
 |f-\ftilde^{(\nu_2)}|_h
 \leq
 \nu_2^{-1}(10n)^3
 \bigl|
  [f,s(h,h')]
 \bigr|_h
 \leq
 \epsilon_1(n,A,\rho)
 (10n)^3
   \leq
   \epsilon_0(n,A,\rho).
\]
By Lemma \ref{lem;23.2.22.20},
there exists an $\ftilde^{(\nu_2)}$-cyclic vector.
But, it contradicts $W^{(\nu_2)}\neq 0$,
according to Corollary \ref{cor;22.9.12.4}.
Hence, we obtain $\max\nbigs(h,h')-1< \nu_2$.
Then, we obtain
$|s-\id|_h\leq n\nu_2\leq C_1(n,A,\rho)\epsilon$.
\hfill\qed

\subsection{Proof of Theorem \ref{thm;23.4.2.11}}

Let $X$, $(\hyperk_{X,n},\theta)$, $C$
and $h_i$ $(i=1,2)$ be as in Theorem \ref{thm;23.4.2.11}.
Let $s$ be the automorphism of
$\hyperk_{X,n}$ determined by $h_2=h_1\cdot s$.
We have
\[
 \sqrt{-1}\Lambda_{g_X}\delbar\del\tr(s)
 \leq
 -\bigl|
 \delbar(s)s^{-1/2}
 \bigr|^2_{h_1,g_X}
-\bigl|
 \bigl[s,\theta\bigr]
 s^{-1/2}
 \bigr|^2_{h_1,g_X}.
\]
By Omori-Yau maximum principle,
there exist $m_0\in\seisuu_{>0}$
and a family of points $p_m\in X$ $(m\geq m_0)$
such that
\[
 \tr(s)(p_m)
 \geq
 \sup\tr(s)-\frac{1}{m},
 \quad\quad
 \sqrt{-1}\Lambda_{g_X}\delbar\del\tr(s)
 \geq-\frac{1}{m}.
\]
Because $h_1$ and $h_2$ are mutually bounded,
there exists $C_1>0$ such that
\[
  \bigl|
 \bigl[s,\theta\bigr]
 \bigr|^2_{h_1,g_X}(p_m)
 \leq
 \frac{C_1}{m}.
\]

Let $\tau_{m}$ be a frame of
the cotangent space of $X$ at $p_m$
such that $|\tau_m|_{g_X}=1$.
It induces a frame
$e_{m,j}=\tau_m^{(n+1-2j)/2}$ $(j=1,\ldots,n)$
of $\hyperk_{X,n|p_m}$.
Because both $h_i$ are mutually bounded with $h_X$,
there exists a constant $B>0$ such that
$|e_{m,1}|_{h_i}\leq B$ for any $m$ and $i$.
Let $f_m$ be the endomorphism of
$\hyperk_{X,n|p_m}$ determined by
$\theta_{|p_m}=f_m\,\tau_m$.
Because $\theta$ is bounded with respect to $h_i$ and $g_X$,
there exists $C_2>0$ independently from $m$
such that $|f_m|_{h_i}\leq C_2$.
By Theorem \ref{thm;22.9.12.10},
there exists $C_3>0$
independently from $m$
such that
\[
 \bigl|
 s-\id
 \bigr|_{h_1}(p_m)
 \leq
 \frac{C_3}{\sqrt{m}}.
\]
Because both $h_i$ are compatible with
the non-degenerate pairing $C$,
we have $\det(s)=1$.
There exists $C_4>0$ independently from $m$
such that
\[
n\leq \sup_X \tr(s)
\leq
 \tr (s)(p_m)
+\frac{1}{m}
\leq n+\frac{C_4}{\sqrt{m}}+\frac{1}{m}.
\]
We obtain that
$\tr(s)$ is constantly $n$,
i.e., $s=\id$.
\hfill\qed

\section{Hitchin section for $SL(n,\mathbb R)$}\label{HitchinSection}

\subsection{Existence of weakly dominant harmonic metric in the general case}

Given a tuple of holomorphic differentials $\vecq=(q_2,q_3,\cdots, q_n)$, one can construct a $SL(n,\mathbb R)$-Higgs bundle 
\[\Big(\hyperk_{X,n}=K_X^{\frac{n-1}{2}}\oplus K_X^{\frac{n-3}{2}}\oplus \cdots\oplus K_X^{\frac{3-n}{2}}\oplus K_X^{\frac{1-n}{2}},\quad \theta(\vecq)=\begin{pmatrix}0&q_2&q_3&q_4&\cdots&q_n\\
1&0&q_2&q_3&\ddots&\vdots\\
&1&0&q_2&\ddots&\vdots\\
&&\ddots&\ddots&\ddots&q_3\\
&&&\ddots&\ddots&q_2\\
&&&&1&0
\end{pmatrix}
\Big).\]
The natural pairings
$K_X^{(n-2i+1)/2}
\otimes
K_X^{-(n-2i+1)/2}
\to \nbigo_X$
induce a non-degenerate symmetric bilinear form
$C_{\hyperk,X,n}$ of
$\hyperk_{X,n}$.
It is a non-degenerate symmetric pairing of
$(\hyperk_{X,n},\theta(\vecq))$. 
Such Higgs bundles are called Higgs bundles in the Hitchin section. They were first introduced by Hitchin in \cite{Hit92} for compact hyperbolic Riemann surfaces. There are various expressions of Higgs bundles in the Hitchin section and they are equivalent to each other. One may refer Appendix \ref{CharacterizationHitchinSection} for details. 

A non-compact Riemann surface $X$ is called hyperbolic if its universal cover is isomorphic to the unit disk $\mathbb D$. A non-compact Riemann surface $X$ is hyperbolic iff it is not $\mathbb C$ nor $\mathbb C^*$. Suppose $X$ is hyperbolic. 

Let $g_X$ be the unique complete conformal hyperbolic metric on $X$. Locally, write $g_X=g_0(dx^2+dy^2)$. The induced Hermitian metric on $K_X^{-1}$ is $\frac{g_0}{2}dz\otimes d\bar z$, also denoted by $g_X$. Denote by $F(g_X)$ the curvature of the Chern connection of the Hermitian metric $g_X$ on $K_X^{-1}$. So $F(g_X)=\bar\partial\partial\log \frac{g_0}{2}=-\partial\bar\partial\log g_0$. The Gaussian curvature of $g_X$ is $k_{g_X}:=\sqrt{-1}\Lambda_{g_X} F(g_X)=-\frac{2}{g_0}\partial_z\partial_{\bar z}\log g_0=-\frac{1}{2}\triangle_{g(X)}\log g_0$. Here $g_X$ is hyperbolic means $k_{g_X}=-1.$

Let $F_i=\oplus_{k=1}^iK_X^{\frac{n+1-2k}{2}}.$ Thus $\mathbf{F}=\{F_1\subset F_2\subset\cdots\subset F_n\}$ forms a full holomorphic filtration of $\hyperk_{X,n}$. And $\theta(\vecq)$ takes $F_i$ to $F_{i+1}\otimes K_X$ and induces an isomorphism between $F_i/F_{i-1}\rightarrow F_{i+1}/F_i\otimes K_X$ for $i=1,\cdots, n-1$. Then, $(\hyperk_{X,n},\theta(\mathbf{0}))$ is the graded Higgs bundle of $(\hyperk_{X,n},\theta(\vecq))$ with respect to the filtration $\mathbf{F}.$

Let $$h_X=\oplus_{k=1}^na_{k,n}\cdot g_X^{-\frac{n+1-2k}{2}},$$ where 
\begin{equation}
\label{a_kn}
a_{k,n}=\prod_{l=1}^{k-1}(\frac{l(n-l)}{2})^{\frac{1}{2}}\cdot \prod_{l=k}^{n-1}(\frac{l(n-l)}{2})^{-\frac{1}{2}}.
\end{equation}
One may check that $h_X$ is a harmonic metric for the Higgs bundle $(\hyperk_{X,n}, \theta(\mathbf{0}))$. 

We call a Hermitian metric $h$ on $\hyperk_{X,n}$ \textbf{weakly dominates} $h_X$ if $\det(h|_{F_k})\leq \det(h_X|_{F_k})$ for $1\leq k\leq n-1.$ 

\begin{thm}\label{existence}
On a hyperbolic surface $X$, there exists a harmonic metric $h$ on $(\hyperk_{X,n},\theta(\vecq))$ satisfying (i) $h$ weakly dominates $h_X$; (ii) $h$ is compatible with $C_{\hyperk,X,n}.$

Moreover, the norm of Higgs field satisfies $|\theta(\vecq)|_{h,g_X}^2\geq |\theta(\mathbf{0})|_{h_X,g_X}^2=\frac{n(n^2-1)}{12}.$

As a result, the associated harmonic map $f: (\widetilde X, \widetilde{g_X})\rightarrow SL(n,\mathbb R)/SO(n)$ satisfies the energy density $e(f)\geq \frac{n^2(n^2-1)}{6}.$ The equality holds if $\vecq=0.$
\end{thm}
\pf
The existence follows from Part (i) of Theorem \ref{MainTheorem1}.

The proof of the moreover statement is identical to the one in \cite[Theorem 4.2]{Li19}. 

From \cite[Section 5.2]{LiIntroduction}, we know that the energy density is $e(f)=2n\cdot |\theta(\vecq)|_{h,g_X}^2$. So $e(f)\geq 2n\cdot |\theta(\mathbf{0})|_{h_X,g_X}^2=\frac{n^2(n^2-1)}{6}.$
\hfill\qed

\subsection{Uniqueness in the case of bounded differentials}
Next, we consider the case when $q_i (i=2,\cdots,n)$ are bounded with respect to $g_X$, that is, $(q_i\bar q_i)/g_X^i$ is bounded.

\begin{thm}\label{BoundedDifferentialExistence}
On a hyperbolic surface $X$, suppose $q_i (i=2,\cdots,n)$ are bounded with respect to $g_X$. Then there uniquely exists a harmonic metric $h$ of $(\hyperk_{X,n},\theta(\vecq))$ over $X$ such that (i) $h$ weakly dominates $h_X$, (ii) $h$ is compatible with $C_{\hyperk, X,n}$.

Moreover, $h$ is mutually bounded with $h_X$.
\end{thm}
\pf
The existence follows from Theorem \ref{existence}.
Let $h_i$ $(i=1,2)$
be harmonic metrics of $(\hyperk_{X,n},\theta(\vecq))$
compatible with $C_{\hyperk,X,n}$
which weakly dominate $h_X$.
By Proposition \ref{prop;23.4.2.10},
both $h_i$ are mutually bounded with $h_X$.
By Theorem \ref{thm;23.4.2.11},
we obtain $h_1=h_2$.
\hfill\qed



\begin{rem}\label{ReplaceCondition}
 The condition (i) in Theorem \ref{BoundedDifferentialExistence} can be replaced by (i') there exists a positive constant $c$
such that
$h_{|F_1(\hyperk_{X,n})}\leq c\cdot h_{X|F_1(\hyperk_{X,n})}$.
\end{rem}

\subsubsection{Compact case}

We reprove the existence and uniqueness of a harmonic metric on $(\hyperk_{X,n},\theta(\vecq))$ over a compact hyperbolic Riemann surface. Note that here our proof does not invoke the Hitchin-Kobayashi correspondence by using the stability of Higgs bundle. 
\begin{thm}\label{CompactExistence}
Given a tuple of holomorphic differentials $\vecq=(q_2,\cdots, q_n)$ on a compact hyperbolic surface $X$, there uniquely exists a harmonic metric $h$ on $(\hyperk_{X,n},\theta(\vecq))$ satisfying $h$ is compatible with $C_{\hyperk,X,n}.$

Moreover, $h$ weakly dominates $h_X$;
\end{thm}
\pf
We first show the existence. Let $X$ be covered by $\mathbb D$ under the map $p:\mathbb D\rightarrow X$, with the covering transformation group of $X$ be $\Gamma<Aut(\mathbb D)=PSL(2,\mathbb R)$, i.e. $X=\mathbb D/\Gamma$. Lift $\vecq, g_X, h_X,\hyperk_{X,n},\theta(\vecq), C_{\hyperk, X,n} $ to $\hat \vecq, g_{\mathbb D}, h_{\mathbb D}, \hyperk_{\mathbb D,n},\theta(\hat\vecq), C_{\hyperk, \mathbb D,n}$ on $\mathbb D$, which are invariant under $\Gamma$. By Theorem \ref{BoundedDifferentialExistence}, there exists a harmonic metric $\hat h\in  \Harm^{dom}(\hyperk_{\mathbb D,n},\theta(\vecq),C_{\hyperk,\mathbb D, n}: h_{\mathbb D})$. From \S\ref{Pullback}, each $\gamma\in \Gamma$ induces an automorphism on $\Harm^{dom}(\hyperk_{\mathbb D,n},\theta(\vecq), C_{\hyperk,\mathbb D, n}: h_{\mathbb D})$. By the uniqueness in Theorem \ref{BoundedDifferentialExistence}, $\gamma^*(\hat h)=\hat h,$ for $\gamma\in \Gamma$. Hence $\hat h$ descends to a harmonic metric $h$ on $(\hyperk_{X,n} ,\theta(\vecq))$ over $\mathbb D/\Gamma=X$.

The lifted $\hat\vecq$ are bounded with respect to $g_{\mathbb D}$ and any lifted harmonic metric $\hat h$ satisfies there exists a positive constant $c$
such that
$\hat h_{|F_1(\hyperk_{\mathbb D,n})}\leq c\cdot  h_{\mathbb D|F_1(\hyperk_{\mathbb D,n})}$. By Theorem \ref{BoundedDifferentialExistence} and Remark \ref{ReplaceCondition}, $\hat h$ is unique and weakly dominates $h_{\mathbb D}$. Thus the descended $h$ is unique and weakly dominates $h_X$. 
\hfill\qed

\subsubsection{Pull back}\label{Pullback}

Let $F:X_1\lrarr X_2$ be a holomorphic map
of Riemann surfaces
which is locally an isomorphism,
i.e., the derivative of $F$ is nowhere vanishing.
Let $\vecq=(q_2,\cdots, q_n)$ be a tuple of holomorphic  differentials on $X_2$.

Because $F$ is locally an isomorphism,
there exists a natural isomorphism
$$F^{\ast}(\hyperk_{X_2,n},\theta(\vecq), C_{\hyperk, X_2,n})
\simeq
(\hyperk_{X_1,n},\theta(F^{\ast}\vecq), C_{\hyperk, X_1,n}).$$
For any harmonic metric $h$ of $(\hyperk_{X_2,n},\theta(F^{\ast}\vecq))$ compatible with $C_{\hyperk, X_2,n}$,
it is well known and easy to check that
the induced metric $F^{\ast}(h)$
of $\hyperk_{X_1,r}$ is a harmonic metric of 
$(\hyperk_{X_1,r},\theta(F^{\ast}\vecq))$ compatible with $C_{\hyperk, X_1,n}$.

Let $h_0$ be a Hermitian metric on $\hyperk_{X_2,n}$. If $h$ weakly dominates $h_0$, then $F^*(h)$ weakly dominates $F^*(h_0)$. Let $h_0$ be a Hermitian metric on $\hyperk_{X,n}$.

In this way,
we obtain the map
\[
F^{\ast}:\Harm^{dom}(\hyperk_{X_2,n},\theta(\vecq), C_{\hyperk, X_2,n}: h_0)\lrarr \Harm^{dom}(\hyperk_{X_1,n},\theta(F^*\vecq), C_{\hyperk, X_1,n}: F^*(h_0)).
\]
If $X_1=X_2$ and $F^{\ast}(\vecq)=\vecq$, $F^*h_0=h_0$,
then $F$ induces an automorphism on $\Harm^{dom}(\hyperk_{X_1,n},\theta(F^*\vecq): h_0)$.

\section{Existence with bounded condition on the unit disk}\label{BoundedExistenceSection}
In this section, let $X$ be the unit disk $\{z\in \mathbb C|~|z|<1\}.$ 
\subsection{Some function spaces}
Let $\mathcal A$ be the set consisting
of all smooth nonnegative functions $f$ such that 
$$\int_X f(z)(1-|z|^2)d\sigma<\infty,$$ where $d\sigma$ is the Lebesgue measure on the unit disk $X$.

Let $G(z,\xi)$ denote the Green function in $X$. Equivalently, from Lemma \ref{PropertiesGreenFunction}, $\mathcal A$ is the set consisting
of all smooth nonnegative functions $f$ such that for some (thus for all) $z$,
$$\int_X G(z,\xi)f(\xi)d\sigma_\xi<\infty.$$

Let $\mathcal A^b$ be the set consisting
of all smooth nonnegative functions $f$ such that 
$$\sup_{z\in X}\int_X G(z,\xi)f(\xi)d\sigma_\xi<\infty.$$ It is clear that $\mathcal A^b\subset \mathcal A.$ From Lemma \ref{PropertiesGreenFunction}, for $p>-2$, $(1-|z|^2)^p\in\mathcal A^b;$ for $p\leq -2,$ $(1-|z|^2)^p\notin\mathcal A.$

\subsection{General existence with bounded condition}
We set $X=\{z\in \mathbb C|~|z|<1\}$ with the Poincar\'e metric $g_X$ and Euclidean metric $g_0(X)$:
$$g_X=\frac{dx^2+dy^2}{(1-|z|^2)^2},\quad g_0(X)=dx^2+dy^2.$$

\begin{prop}\label{BasicExistence} Suppose $(E,\delbar_E=\delbar_E^0+\xi, \theta=\theta_0+\phi)$ is a Higgs bundle over $X$ for $\phi\in A^{1,0}(X, \End(E))$ and $\xi\in A^{0,1}(X, \End(E))$. Assume $(E,\delbar_E^0, \theta_0)$ is a Higgs bundle over $X$ which admits a harmonic metric $h_1$. 
\begin{itemize}
\item Suppose $$|[\phi,(\theta_0)^{*h_1}]|_{h_1, g_0(X)}\in \mathcal A,\quad  |\phi|_{h_1, g_0(X)}^2\in \mathcal A,\quad |\delbar_E^{0}\xi^{*h_1}|_{h_1,g_0(X)}\in \mathcal A, \quad |\xi|_{h_1, g_0(X)}^2\in \mathcal A.$$
Then there exists a harmonic metric $h$ on $(E,\delbar_E, \theta)$.
\item 
Suppose  $$|[\phi,(\theta_0)^{*h_1}]|_{h_1, g_0(X)}\in \mathcal A^b,\quad  |\phi|_{h_1, g_0(X)}^2\in \mathcal A^b,\quad |\delbar_E^{0}\xi^{*h_1}|_{h_1,g_0(X)}\in \mathcal A^b, \quad |\xi|_{h_1, g_0(X)}^2\in \mathcal A^b.$$
Then there exists a harmonic metric $h$ on $(E,\delbar_E, \theta)$ mutually bounded with $h_1$.
\end{itemize}
\end{prop}
\pf
Let $\nabla_{h_1}=\del_E^{h_1}+\delbar_E^0$ be the Chern connection of $E$ determined by $h_1$ and $\delbar_E^0$. We have 
\begin{eqnarray*}
F(\delbar_E^0+\xi, \theta_0+\phi,h_1)&=&F(\nabla_{h_1})+\del_E^{h_1}\xi-\delbar_E^0\xi^{*h_1}-[\xi, \xi^{*h_1}]+[\theta_0+\phi, (\theta_0+\phi)^{*h_1}]\\
&=&-[\theta_0, {\theta_0}^{*h_1}]+\del_E^{h_1}\xi-\delbar_E^0\xi^{*h_1}-[\xi, \xi^{*h_1}]+[\theta_0+\phi, (\theta_0+\phi)^{*h_1}]\\
&=&\del_E^{h_1}\xi-\delbar_E^0\xi^{*h_1}-[\xi, \xi^{*h_1}]-[\phi, {\theta_0}^{*h_1}]-[\theta_0, \phi^{*h_1}]-[\phi, \phi^{*h_1}].
\end{eqnarray*}

By assumption,
\begin{eqnarray*}
\Big|\Lambda_{g_0(X)}F(\delbar_E^0+\xi, \theta_0+\phi,h_1)\Big|_{h_1}\in \mathcal A.
\end{eqnarray*}

For $0<r<1,$ we set $X_r:=\{|z|<r\}.$ Let $h_{X_r}$ be the harmonic metric of $(E,\theta)|_{X_r}$ such that $h_{X_r}=h_1$ on $\partial X_r$. We have $\det(s(h_1|_{X_r},h_{X_r}))=1.$ Recall $\triangle=\partial_x^2+\partial_y^2=4\partial_z\partial_{\bar z}$. We have the following inequality on $X_r:$
$$\frac{1}{2}\triangle \log \tr(s(h_1|_{X_r},h_{X_r}))= \sqrt{-1}\Lambda_{g_0(X)}\partial\bar\partial\log \tr(s(h_1|_{X_r},h_{X_r}))\geq -\Big|\Lambda_{g_0(X)}F(\theta+\phi,h_1)\Big|_{h_1}.$$

\begin{lem} \label{DirichletProblem}
Let $f$ be a nonnegative smooth function on $X$.
Suppose $f\in \mathcal A.$
Let $u_r$ be the unique solution satisfying 
\begin{eqnarray*}
    &&\triangle u_r=-f,\quad \text{in $X_r$}\\
    &&u_r=0\quad \text{on $\partial X_r$}
\end{eqnarray*}
There exists a smooth nonnegative function $v$ on $X$ such that $0\leq u_r\leq v,$ for any $r\in (0,1)$.

If moreover $f\in \mathcal A^b$, we can choose a bounded $v$ satisfying the above property.
\end{lem}
\pf
Define $v(z)=\frac{1}{2\pi}\int_{X}f(\xi)G(z,\xi)d\sigma_\xi.$ By Lemma \ref{PropertiesGreenFunction}, we know $v(z)$ is well-defined and nonnegative. 
Then we have \begin{eqnarray*}
    &&\triangle v=-f,\quad \text{in $X$}\\
    &&v\geq 0\quad \text{in $X$}
\end{eqnarray*}
By the maximum principle, $u_r\leq v$ holds in $X_r$.
Also, note that $0$ is a subsolution to the equation. By the maximum principle, $u_r\geq 0.$

If $f\in \mathcal A^b,$ then $v$ is bounded by definition. 
\hfill\qed\\

By Lemma \ref{DirichletProblem} and using the maximum principle, there exists a smooth function $v$ on $X$ such that $$\log (\tr(s(h_1|_{X_r},h_{X_r}))/\text{rank}(E))\leq v.$$

Then, by Proposition \ref{prop;20.5.29.1}, there is a convergence subsequence of $h_{X_r}$ whose limit is denoted by $h$. Then $h$ is a harmonic metric of $(E,\theta).$  

If $v$ is bounded, $h$ is mutually bounded with $h_1$.
\hfill\qed

\subsection{Existence for holomorphic chains}
Given $k$ holomorphic vector bundles $E_i$ over $X$ of rank $n_i$, $i=1,\cdots,k$. We can consider a Higgs bundle as following:
 $$(E=E_1\oplus E_2\oplus\cdots\oplus E_k,\quad \theta=\begin{pmatrix}0&&&&\\\theta_1&0&&&\\&\theta_2&0&&\\&&\ddots&\ddots&\\&&&\theta_{k-1}&0\end{pmatrix}),$$ where $\theta_i\in H^0(X,\Hom(E_i,E_{i+1})\otimes K).$ Such Higgs bundle is called a holomorphic chain of type $(n_1,\cdots, n_k)$. We call a Hermitian metric $h$ on $E$ orthogonal if $h$ is orthogonal with respect to the decomposition $E=\oplus_{i=1}^k E_i$. 
 
For such $(E, \theta)$, we also consider a holomorphic chain $(E,\theta_0)$ as follows: $$(E,\quad  \theta_0=\begin{pmatrix}0&&&&\\\phi_1&0&&&\\&\phi_2&0&&\\&&\ddots&\ddots&\\&&&\phi_{k-1}&0\end{pmatrix}),$$ where there exists a subset $S$ of $\{1,2,\cdots, k-1\}$, $\phi_i=0$ for $i\in S$ and $\phi_i=\theta_i$ for $i\notin S$.

In the following, we will deduce the existence of a harmonic metric on $(E,\theta)$ from the one on $(E,\theta_0)$ if it exists. 
\begin{thm}\label{NilpotentHiggsBundles}
We consider two holomorphic chains $(E,\theta), (E, \theta_0)$ as above. Suppose there exists an orthogonal harmonic metric $h_1$ on $(E, \theta_0)$.
Suppose $|\theta_i|_{h_1, g_0(X)}^2\in \mathcal A (i\in S),$ then there exists an orthogonal harmonic metric $h$ on $(E, \theta)$.

Moreover, $|\theta_i|_{h_1, g_0(X)}^2\in \mathcal A^b(i\in S)$ if and only if there exists an orthogonal harmonic metric $h$ on $(E,\theta)$ mutually bounded with $h_1$.
\end{thm}
\pf
Note that $[\theta-\theta_0, (\theta_0)^{*h_1}]=0.$ The existence part under both assumptions follows from Proposition \ref{BasicExistence}.

We only need to prove the inverse direction for a bounded metric. Suppose we have a diagonal harmonic metric $h$ mutually bounded with $h_1$.

The Hitchin equation for $h_1$ on $(E_1\oplus E_2\oplus \cdots \oplus E_k, \theta_0)$ is 
\begin{eqnarray*}
&&F(h_1|_{E_1})=- (\phi_1)^{*h_1}\wedge\phi_1\\
&&F(h_1|_{E_2})=-(\phi_2)^{*h_1}\wedge\phi_2- \phi_1\wedge(\phi_1)^{*h_1}\\
&&\cdots\\
&&F(h_1|_{E_{k-1}})=-(\phi_{k-1})^{*h_1}\wedge\phi_{k-1}-\phi_{k-2}\wedge(\phi_{k-2})^{*h_1}\\
&&F(h_1|_{E_k})=-\phi_{k-1}\wedge(\phi_{k-1})^{*h_1}
\end{eqnarray*}
Denote by $$|\phi_i|_{h_1,g_0(X)}^2=\sqrt{-1}\Lambda_{g_0(X)}\tr(\phi_i\wedge(\phi_i)^{*h_1})=-\sqrt{-1}\Lambda_{g_0(X)}\tr((\phi_i)^{*h_1}\wedge\phi_i).$$ So
\begin{eqnarray*}
&&-\sqrt{-1}\Lambda_{g_0(X)}\tr (F(h_1|_{E_1}))=-|\phi_1|_{h_1,g_0(X)}^2\\
&&-\sqrt{-1}\Lambda_{g_0(X)}\tr (F(h_1|_{E_2}))=-|\phi_2|_{h_1,g_0(X)}^2+|\phi_1|_{h_1,g_0(X)}^2\\
&&\cdots\\
&&-\sqrt{-1}\Lambda_{g_0(X)}\tr(F(h_1|_{E_{k-1}}))=-|\phi_{k-1}|_{h_1,g_0(X)}^2+|\phi_{k-2}|_{h_1,g_0(X)}^2\\
&&-\sqrt{-1}\Lambda_{g_0(X)}\tr(F(h_1|_{E_k}))=|\phi_{k-1}|_{h_1,g_0(X)}^2
\end{eqnarray*}

Therefore, \begin{equation*}
-\sqrt{-1}\Lambda_{g_0(X)}\tr (F(h_1|_{\oplus_{l=1}^iE_l}))=-|\phi_i|_{h_1,g_0(X)}^2,\quad i=1,2,\cdots, k-1.\end{equation*}

Similarly, 
\begin{equation*}
-\sqrt{-1}\Lambda_{g_0(X)}\tr (F(h|_{\oplus_{l=1}^iE_l}))=-|\theta_i|_{h,g_0(X)}^2,\quad i=1,2,\cdots, k-1.\end{equation*}

Let $u_i=\log(\det (h|_{\oplus_{l=1}^iE_l})/\det(h_1|_{\oplus_{l=1}^iE_l})).$
Then
$$-\sqrt{-1}\Lambda_{g_0(X)}\tr(F(h|_{\oplus_{l=1}^iE_l}))+\sqrt{-1}\Lambda_{g_0(X)}\tr(F(h_1|_{\oplus_{l=1}^iE_l}))=2\partial_z\partial_{\bar z}u_i=\frac{1}{2}\triangle u_i.$$

Then we obtain the equation for $u_i$'s as follows:
\begin{equation*}
\frac{1}{2}\triangle u_i=-|\theta_i|_{h,g_0(X)}^2+|\phi_i|_{h_1,g_0(X)}^2,,\quad i=1,2,\cdots, k-1.\end{equation*}

 We only focus on $u_i$'s where $i\in S$. For $i\in S$, $\phi_i=0$ and the equation becomes:
\begin{equation*}
\frac{1}{2}\triangle u_i=-|\theta_{i}|_{h,g_0(X)}^2
\end{equation*}
Suppose $|u_i|\leq M$ for  $i=1,2,\cdots, k$. 
Applying Lemma \ref{ConverseLemma}, we obtain $|\theta_{i}|_{h,g}^2\in \mathcal A^b$ for $i\in S.$ Since $|u_i|\leq M$ for  $i=1,2,\cdots, k$, $|\theta_{i}|_{h_1,g}^2\in \mathcal A^b$ for $i\in S.$
\hfill\qed\\

In the following, we will show Lemma \ref{ConverseLemma} which was used in proving Theorem \ref{NilpotentHiggsBundles}.

Let $u$ be a subharmonic function on a domain $D$. A harmonic majorant of $u$ is a harmonic function $h$ on $D$ such that $h\geq u$ there. If also $h\leq k$, for every other harmonic majorant $k$ of $u,$ then $h$ is called the least harmonic majorant of $u$.
\vspace{0.1cm}

\begin{lem}(\cite[Theorem 3.3]{TopicsMajorant})\label{HarmonicMajorant}
Let $u$ be a subharmonic function on $X$ with $u\neq -\infty.$ There exists a harmonic majorant for $u$ if and only if 
$$\sup\limits_{0<r<1}\frac{1}{2\pi}\int_{0}^{2\pi}u(re^{it})dt<\infty.$$
\end{lem}

\begin{lem}{\cite[Theorem 4.5.4]{Ransford}}\label{PoissonJensen}
Let $u$ be a subharmonic function on $X$ such that $u\neq -\infty.$ If $u$ has a harmonic majorant on $X$, then it has a least one, $h$, and 
$$u(z)=h(z)-\frac{1}{2\pi}\int_XG(z,\xi)\triangle u(\xi)d\sigma_\xi.$$
\end{lem}
\begin{lem}\label{ConverseLemma} Let $f$ be a nonnegative or nonpositive smooth function.
If $|u|\leq M$ and $\triangle u=f$ on $X$, then $|f|\in \mathcal A^b.$
\end{lem}
\pf
It is enough to show the case for $f$ being nonnegative. If $f$ is nonpositive, we can consider $\triangle (-u)=-f.$

By Lemma \ref{HarmonicMajorant}, since $|u|\leq M$, there exists a least harmonic majorant $h$ for $u$ and $h\leq M$ since the constant function $M$ is a harmonic majorant of $u$.

By Lemma \ref{PoissonJensen}, 
\begin{equation*}
\frac{1}{2\pi}\int_{X} G(z,\xi)\triangle u(\xi)d\sigma_{\xi}=h(z)-u(z)\leq 2M.
\end{equation*}
Thus $f=\triangle u\in \mathcal A^b.$
\hfill\qed

\subsection{Relation to prescribed curvature equation}
Consider the curvature equation on $X$: 
\begin{equation}\label{CurvatureEquation}
\frac{1}{4}\triangle u=|\alpha|^2e^{2u},
\end{equation}where $\alpha$ is a holomorphic function on $X$. That is, we are looking for the function $u$ such that the metric $e^{2u}(dx^2+dy^2)$ on $X$ has Gaussian curvature $-4|\alpha|^2.$

As a corollary of Theorem \ref{NilpotentHiggsBundles}, we can recover the following theorem shown by Kraus. 
\begin{prop}(\cite[Theorem 3.1]{Kraus})\label{boundedRank2}
(1) $|\alpha|^2\in \mathcal A^b$ if and only if there exists a solution $u$ bounded from above and below, of Equation (\ref{CurvatureEquation}).

(2) if $|\alpha|^2\in \mathcal A$, then there exists a solution $u$ of Equation (\ref{CurvatureEquation}).
\end{prop}
\pf
Consider the Higgs bundle $(\mathcal O\oplus \mathcal O, \begin{pmatrix}0&0\\\alpha &0\end{pmatrix}dz)$ over the unit disk $X$. Note that the Higgs bundle is symmetric to the non-degenerate symmetric pairing $C=\begin{pmatrix}0&1\\1&0\end{pmatrix}$. On $\mathcal O\oplus\mathcal O$, there is a flat Hermitian metric $h_1=\diag(1,1)$ which is symmetric with respect to $C$ and of unit determinant. Then we obtain $|\alpha\cdot dz|_{h_1,g_0(X)}^2=2|\alpha|^2.$ Also, a diagonal harmonic metric which is compatible with $C$ is of unit determinant. It is of the form $h=\diag((h^0)^{-1}, h^0)$. Let $h^0=e^u,$ then $u$ satisfies Equation (\ref{CurvatureEquation}).

The rest follows from Theorem \ref{NilpotentHiggsBundles}.
\hfill\qed\\

In fact, Kraus showed the converse direction for the existence of the curvature equation. 
\begin{prop}(\cite[ Theorem 1.3]{Kraus})\label{boundedRank2}
If there exists a solution of Equation (\ref{CurvatureEquation}), then there exists a non-vanishing holomorphic function $f$ on $X$ such that $|\alpha\cdot f|^2\in \mathcal A$.
\end{prop}
\begin{rem}
Kraus' proof relies on the Littlewood-Paley identity for holomorphic functions.
It is not clear if such conditions are still necessary for higher rank nilpotent Higgs bundles.
\end{rem}

\subsection{Holomorphic chains of type $(1,1,\cdots,1)$}
The following proposition indicates that only a proper subset of $\theta_i$'s being nice does not imply the existence of a diagonal harmonic metric.
\begin{prop}\label{SomeNoneExistence}
Consider a Higgs bundle $(\mathcal O\oplus\mathcal O\oplus \cdots\oplus \mathcal O, \begin{pmatrix}0&&&&\\\gamma_1&0&&&\\&\gamma_2&0&&\\&&\ddots&\ddots&\\&&&\gamma_{n-1}&0\end{pmatrix}dz)$ over the unit disk $X$ satisfying $\prod_{i=1}^{n-1}\gamma_i^{i(n-i)}=\alpha^{\frac{n(n^2-1)}{6}}$ for a holomorphic function $\alpha$ which is not constantly $0$. A necessary condition for the existence of a diagonal harmonic metric is that there exists a non-vanishing holomorphic function $f$ on $X$ such that $|\alpha\cdot f|^2\in \mathcal A.$
\end{prop}
\pf
Suppose there is a harmonic metric $h=\diag(e^{-u_1},e^{-u_2},\cdots,e^{-u_n})$. Let $w_k=\sum_{i=1}^ku_i.$
Then following from the calculations in the proof of Theorem \ref{NilpotentHiggsBundles} and $|dz|_{g_0(X)}^2=2$, the Hitchin equation becomes
\begin{eqnarray*}
&&\frac{1}{4}\triangle w_1=|\gamma_1|^2e^{2w_1-w_2}\\
&&\frac{1}{4}\triangle w_2=|\gamma_2|^2e^{2w_2-w_1-w_3}\\
&&\cdots\\
&&\frac{1}{4}\triangle w_{n-2}=|\gamma_{n-2}|^2e^{2w_{n-2}-w_{n-3}-w_{n-1}}\\
&&\frac{1}{4}\triangle w_{n-1}=|\gamma_{n-1}|^2e^{2w_{n-1}-w_{n-2}}
\end{eqnarray*}
Summing up the above $(n-1)$-equations, we obtain
\begin{equation}\label{SumTogether}
\frac{1}{4}\triangle(w_1+w_2+\cdots+w_{n-1})=|\gamma_1|^2e^{2w_1-w_2}+\sum_{i=2}^{n-2}|\gamma_i|^2e^{2w_i-w_{i-1}-w_{i+1}}+|\gamma_{n-1}|^2e^{2w_{n-1}-w_{n-2}}
\end{equation}

Let $r_i=\frac{i(n-i)}{2}.$ Then $$2r_1-r_2=1, 2r_i-r_{i-1}-r_{i+1}=1(i=2,\cdots,n-2),2r_{n-1}-r_{n-2}=1, \sum_{i=1}^{n-1}r_i=\frac{n(n^2-1)}{12}.$$

Note that the right hand side of Equation (\ref{SumTogether}) satisfies 
\begin{eqnarray*}
&&|\gamma_1|^2e^{2w_1-w_2}+\sum_{i=2}^{n-2}|\gamma_i|^2e^{2w_i-w_{i-1}-w_{i+1}}+|\gamma_{n-1}|^2e^{2w_{n-1}-w_{n-2}}\\
\geq&&\frac{1}{\max_{i=1,\cdots,n-1}r_i}\cdot(r_1|\gamma_1|^2e^{2w_1-w_2}+\sum_{i=2}^{n-2}r_i|\gamma_i|^2e^{2w_i-w_{i-1}-w_{i+1}}+r_{n-1}|\gamma_{n-1}|^2e^{2w_{n-1}-w_{n-2}})\\
\geq&&\frac{1}{\max_{i=1,\cdots,n-1}r_i}\cdot \sum_{i=1}^{n-1}r_i\cdot\big(\prod_{i=1}^{n-1}|\gamma_i|^{2r_i}e^{r_1(2w_1-w_2)+\sum_{i=2}^{n-2}r_i(2w_i-w_{i-1}-w_{i+1})+r_{n-1}(2w_{n-1}-w_{n-2})}\big)^{\frac{1}{\sum_{i=1}^{n-1}r_i}}\\
\geq&&\frac{1}{\max_{i=1,\cdots,n-1}r_i}\cdot\frac{n(n^2-1)}{12}\cdot\big(\prod_{i=1}^{n-1}|\gamma_i|^{2r_i}e^{(2r_1-r_2)w_1+\sum_{i=2}^{n-2}(2r_i-r_{i-1}-r_{i+1})w_i+(2r_{n-1}-r_{n-2})w_{n-1}}\big)^{\frac{12}{n(n^2-1)}}\\
\geq&&\frac{1}{\max_{i=1,\cdots,n-1}r_i}\cdot\frac{n(n^2-1)}{12}\cdot\big(\prod_{i=1}^{n-1}|\gamma_i|^{2r_i}e^{w_1+\cdots+w_{n-1}}\big)^{\frac{12}{n(n^2-1)}}.
\end{eqnarray*}

So 
\begin{equation}
\frac{1}{4}\triangle (w_1+\cdots+w_{n-1})\geq \frac{1}{\max_{i=1,\cdots,n-1}r_i}\cdot\frac{n(n^2-1)}{12} \big(\prod_{i=1}^{n-1}|\gamma_i|^{i(n-i)}\big)^{\frac{12}{n(n^2-1)}}\cdot e^{\frac{12}{n(n^2-1)}(w_1+\cdots+w_{n-1})}. 
\end{equation}

Consider the equation 
\begin{equation}\label{TestEquation}
\frac{1}{4}\triangle u= \frac{1}{\max_{i=1,\cdots,n-1}r_i}\big(\prod_{i=1}^{n-1}|\gamma_i|^{i(n-i)}\big)^{\frac{12}{n(n^2-1)}}\cdot e^{u}. 
\end{equation}
Then $\frac{12}{n(n^2-1)}\cdot(w_1+\cdots+w_{n-1})$ is a subsolution to the equation (\ref{TestEquation}). Note that $\gamma_i$'s are holomorphic functions and not constantly zero. So the function $-\big(\prod_{i=1}^{n-1}|\gamma_i|^{i(n-i)}\big)^{\frac{12}{n(n^2-1)}}$ satisfies the essential negative property in \cite[Definition 0.1]{KalkaYang}.
By \cite[Theorem 4]{KalkaYang}, the existence of a subsolution implies there exists a $C^2$ solution to the equation (\ref{TestEquation}).

The rest follows from Proposition \ref{boundedRank2} and the assumption $\prod_{i=1}^{n-1}\gamma_i^{i(n-i)}=\alpha^{\frac{n(n^2-1)}{6}}$ for a holomorphic function $\alpha$. 
\hfill\qed

\begin{rem}
It would be interesting if one could find a necessary and sufficient condition on $\gamma_i (i=1,\cdots,n-1)$ for the existence of a diagonal harmonic metric.
\end{rem}

\section{$SO(n,n+1)$-Higgs bundles}\label{SO(n,n+1)}
In this section, we discuss the existence of harmonic metrics on $SO(n,n+1)$-Higgs bundles over non-compact Riemann surfaces
by using the techniques developed in this paper and
our previous paper \cite{LiMochizukiGeneric}.

\begin{df}\label{SOBundle}
\begin{itemize}
 \item An $SO(n,n+1)$-Higgs bundle over a Riemann surface $X$ is given by
       $((V,Q_V),(W,Q_W), \eta)$,
       where $(V,Q_V)$ is an orthogonal bundle of rank $n$
       satisfying $\det V=\mathcal O_X$,
       $(W,Q_W)$ is an orthogonal bundle of rank $n+1$
       satisfying $\det W=\mathcal O_X$,
       and $\eta:V\rightarrow W\otimes K_X$ is a holomorphic bundle map. 
 \item The associated $SL(2n+1,\mathbb C)-$Higgs bundle is $$(E,\theta)=\left(V\oplus W, \begin{pmatrix}0&\eta^{\dagger}\\\eta&0\end{pmatrix}\right),$$ where $\eta^{\dagger}$ is the adjoint of $\eta$ with respect to $Q_V, Q_W$. 
 \item A harmonic metric $h$ on $(E,\theta)$ is called compatible with $SO(n,n+1)$-structure if $h=h|_V\oplus h|_W$ where $h_V, h_W$ are compatible with $Q_V, Q_W$ respectively. 
\end{itemize}
\end{df}

\subsection{Dirichlet problem}

Let $Y\subset X$ be a relatively compact connected open subset
with smooth boundary $\del Y$.
Assume that $\del Y$ is non-empty.
Let $h_{\del Y}$ be any Hermitian metric of
$E_{|\del Y}$.
\begin{lem}\label{Symmetry1}
Let $h$ be a harmonic metric of $(E,\theta)$ such that $h|_{\del Y}=h_{\del Y}$. Suppose $h_{\del Y}$ is compatible with $SO(n.n+1)$-structure. Then $h$ is compatible with $SO(n.n+1)$-structure.
\end{lem}
\pf 
First we show that $h=\diag(h_V, h_W).$ There exists the automorphism $\varphi=1_V\oplus (-1_W)$ on $E=V\oplus W$. Because $\varphi^*\theta=-\theta$, $\varphi^*(h)$ is also a harmonic metric of $(E,\theta)$. Because $\varphi^*(h)|_{\partial Y}=h_{\partial Y}$, we obtain $\varphi^*(h)=h$. It means that $h$ is the direct sum of the Hermitian metrics of $V$ and $W$. 

Next we show that $h_V, h_W$ are compatible with $Q_V, Q_W$ respectively. The metric $h$ induces a harmonic metric $h^{\lor}=(h|_V)^{\lor}\oplus (h|_W)^{\lor}$ on $(E^{\lor},\theta^{\lor})$. Let $\Psi_{Q_V}: V\rightarrow V^{\lor}$ and $\Psi_{Q_W}: W\rightarrow W^{\lor}$ be the induced isomorphism by $Q_V, Q_W$, respectively. Then $(\Psi_{Q_V})^*((h|_V)^{\lor})\oplus(\Psi_{Q_W})^*( (h|_W)^{\lor})$ is again a harmonic metric on $(E,\theta)$. Since $\big((\Psi_{Q_V})^*((h|_V)^{\lor})\oplus(\Psi_{Q_W})^*( (h|_W)^{\lor})\big)|_{\del Y}=h_{\del Y}$, we obtain $(\Psi_{Q_V})^*((h|_V)^{\lor})\oplus(\Psi_{Q_W})^*( (h|_W)^{\lor})=h|_V\oplus h|_W$. It means $h_V, h_W$ are compatible with $Q_V, Q_W$ respectively. 
\hfill\qed

\subsection{The generically regular semisimple case}

Let $((V,Q_V),(W,Q_W),\eta)$ be an $SO(n,n+1)$-Higgs bundle on $X$.
Let $(E,\theta)$ be the associated $SL(2n+1,\mathbb C)$-Higgs bundle. We obtain $\eta^{\dagger}\circ \eta\in \End(V)\otimes K_X^2.$ 
Let $(T^*X)^{\otimes 2}$ denote the total space of the line bundle $K_X^2$. Let $Z_X\subset (T^*X)^{\otimes 2}$ denote the zero-section. The spectral curve $\Sigma_{\eta^{\dagger}\circ\eta}\subset (T^*X)^{\otimes 2}$ of $\eta^{\dagger}\circ \eta$ is defined as usual. We obtain the finite map $\pi:\Sigma_{\eta^{\dagger}\circ \eta}\cup Z_X\rightarrow X$. 
\begin{df}
We say that the tuple
$((V,Q_V),(W,Q_W),\eta)$ is generically regular semisimple
if there exists
$P\in X$ such that $|\pi^{-1}(P)|=n+1$.
\end{df}

\begin{thm}\label{GenericExistence}
\label{thm;23.6.16.1}
If $((V,Q_V),(W,Q_W),\eta)$ is generically regular semisimple, then there exists a harmonic metric $h$ of $(E,\theta)$ compatible with $SO(n,n+1)$-structure. 
\end{thm}
\pf
The following lemma follows from Corollary \ref{Equivalence} below.
\begin{lem}\label{TwoEquivalence}
 $((V,Q_V),(W,Q_W),\eta)$ is generically regular semisimple if and only if the associated $SL(2n+1,\mathbb C)$-Higgs bundle $(E,\theta)$ is generically regular semisimple.
(See \cite[Definition]{LiMochizukiGeneric}
for generically regular semisimplicity for Higgs bundles.)
\hfill\qed
\end{lem}

Let $h_0=h_0|_V\oplus h_0|_W$ be a Hermitian metric of $E=V\oplus W$ such that $h_0|_V$ and $h_0|_W$ are compatible with $Q_V$ and $Q_W$, respectively.
Let $X_i$ $(i=1,\cdots)$ be an exhaustion family of $X$. Let $E_i,V_i$ and $W_i$ denote the restriction of $E, V$ and $W$ to $X_i$, respectively.
Let $h_{0,i}$ denote the restriction of $h_0$ to $X_i$.

Let $h_i$ be a harmonic metric of $(E_i,\theta_i)$
such that $h_i|_{\del X_i}=h_0|_{\del X_i}$.
By Lemma \ref{Symmetry1}, $h_i$ is compatible with $SO(n,n+1)$-structure.
It implies that $h_i$ is compatible with
the non-degenerate symmetric pairing $Q_V\oplus Q_W$.
Let $s_{i}$ be the automorphism of $E_i$
determined by $h_i=h_{0,i}\cdot s_i$
as in \S\ref{subsection;23.6.16.10}.
Note that the Higgs field $\theta$ is self-adjoint
with respect to the non-degenerate symmetric pairing
$Q_V\oplus Q_W$ of $E$.
By Lemma \ref{TwoEquivalence}
and \cite[Proposition 2.37]{LiMochizukiGeneric},
there exist positive constants $C_i$ $(i=1,2,\ldots)$
such that the following holds on $X_i$
for $j\geq i+1$:
\[
 \bigl|
 s_j
 \bigr|_{h_{0,i}}
+\bigl|
 s_j^{-1}
 \bigr|_{h_{0,i}}
 \leq C_i.
\]
By Proposition \ref{prop;20.5.29.1},
there exists a convergent subsequence $h_i'$.
As the limit,
we obtain a harmonic metric $h$ of $(E,\theta)$
compatible with $SO(n,n+1)$-structure. 
\hfill\qed

\subsubsection{Appendix: Preliminary from linear algebra}

Let $R$ be any field.
In this subsection,
we consider matrices whose entries are contained in $R$.
For any positive integer $n$,
let $I_n$ denote the $(n\times n)$-identity matrix,
and let $0_n$ denote the $(n\times n)$-zero matrix.

Let $n\geq m$ be positive integers.
Let $A$ be an $(n\times m)$-matrix.
Let $B$ be an $(m\times n)$-matrix.
Let $C$ be the $(n+m)$-square matrix
given as follows:
\[
 C=
 \begin{pmatrix}
  0_n & A \\ B & 0_m
 \end{pmatrix}.
\]
\begin{lem}
We have $\det(tI_{n+m}-C)=t^{n-m}\det(t^2I_m-BA)$
in $R[t]$.
\end{lem}
\pf
It is enough to prove the equality in $R[t,t^{-1}]$.
Let $0_{m,n}$ denote the $(m\times n)$-zero matrix.
We have
\begin{multline}
\det(tI_{n+m}-C)=
 \det
 \begin{pmatrix}
  t I_n & -A \\
  -B & tI_m
 \end{pmatrix}
= \det
 \begin{pmatrix}
  t I_n & -A \\
  0_{m,n} & tI_m-t^{-1}BA
 \end{pmatrix}
\\
=t^{n}\det(tI_m-t^{-1}BA)
=t^{n-m}\det(t^2I_m-BA).
\end{multline}

\hfill\qed

\vspace{.1in}
We recall that an $(\ell\times\ell)$-matrix is called
regular semisimple
if it has $\ell$-distinct eigen values.
\begin{cor}
\label{Equivalence}
If $n\geq m+2$,
$C$ cannot be regular semisimple.
If $n=m,m+1$,
$C$ is regular semisimple
if and only if
$BA$ is invertible and regular semisimple. 
\hfill\qed
\end{cor}

\subsection{Collier section}

Given a holomorphic line bundle $M$ on $X$,
$\mu\in H^0(X, M^{-1}\otimes K_X^n)$,
$\nu\in H^0(X,M\otimes K_X^n)$,
$q_{2i}\in H^0(X, K_X^{2i})$
$(i=1, \cdots, n-1)$,
one can construct the following $SO(n,n+1)$-Higgs bundle
$((V,Q_V),(W,Q_W),\eta_{\mu,\nu}(\vecq))$
given by 
\begin{eqnarray}\label{OriginalBundle1}
&&(V,Q_V)=(K_X^{n-1}\oplus K_X^{n-3}\oplus\cdots\oplus K_X^{3-n}\oplus K_X^{1-n}, \begin{pmatrix}&&1\\&\iddots&\\1&&\end{pmatrix})\nonumber\\
&&(W,Q_W)=(M\oplus K_X^{n-2}\oplus K_X^{n-4}\oplus\cdots\oplus K_X^{4-n}\oplus K_X^{2-n}\oplus M^{-1}, \begin{pmatrix}&&1\\&\iddots&\\1&&\end{pmatrix})\nonumber\\
&&\eta_{\mu,\nu}(\vecq)=\begin{pmatrix}
0&0&0&\cdots&\cdots&0&\nu\\
1&q_2&q_4&\cdots&\cdots&q_{2n-4}&q_{2n-2}\\
&1&q_2&q_4&\cdots&\cdots&q_{2n-4}\\
&&1&q_2&\ddots&\ddots&q_{2n-6}\\
&&&\ddots&\ddots&\ddots&\vdots\\
&&&&\ddots&\ddots&\vdots\\
&&&&&1&q_2\\
&&&&&&\mu
\end{pmatrix}:V\rightarrow W\otimes K_X.
\end{eqnarray}

When $X$ is a compact Riemann surface of genus at lease two, for each integer $d\in (0,n(2g-2)]$, Brian Collier in \cite[Theorem 4.11]{Collier} defined a component $X_d$ of the moduli space of $SO(n,n+1)-$Higgs bundles formed by the above Higgs bundles determined by $(M, \mu, \nu, q_2,\cdots, q_{2n-2})$ where $\deg(M)=d$ and $\mu\neq 0$. In particular, when $d=n(2g-2),$ $X_d$ coincides with the Hitchin component for $SO(n,n+1).$ Such components are analogues of Hitchin components. Such Higgs bundles correspond to positive $SO(n,n+1)$ representations. 

We call the above Higgs bundles are in the Collier section. We are going to discuss the existence of harmonic metrics of Higgs bundles in the Collier section over non-compact Riemann surfaces.

\subsubsection{Existence for the case $\mu\neq 0$}

Let $(E,\theta)$ be the Higgs bundle
associated with
the $SO(n,n+1)$-Higgs bundle $((V,Q_V),(W,Q_W),\eta_{\mu,\nu}(\vecq))$
in (\ref{OriginalBundle1}).
We introduce a holomorphic full filtration
$\mathbf{F}(E)=\{F_1(E)\subset F_2(E)\subset\cdots
\subset F_{2n+1}(E)\}$ of $E$  as follows.
We define
$F_{2i+1}(W)$ $(i=0,\ldots,n)$ by
\[
 F_1(W)=M,
 \quad
 F_{2i+1}(W)=M\oplus
 K_X^{n-2}\oplus\cdots
 \oplus K_X^{n-2j}\,\,\,(i=1,\ldots,n-1),
 \quad
 F_{2n+1}(W)=W,
\]
We also set
$F_{2i}(W)=F_{2i-1}(W)$ for $i=1,\ldots,n$ and $F_0(W)=0$.
We define
$F_{2i}(V)$ $(i=1,\ldots,n)$ by
\[
 F_{2i}(V)=K_X^{n-1}\oplus\cdots
 \oplus K_X^{n+1-2i}\,\,\,(i=1,\ldots,n).
\]
We also set
$F_{2i+1}(V)=F_{2i}(V)$ for $i=1,\ldots,n$
and $F_1(V)=0$.
Then, $\theta$ takes $F_j(W)$ to $F_{j+1}(V)\otimes K_X$,
and $F_{j}(V)$ to $F_{j+1}(W)\otimes K_X$.
We define
\[
 F_j(E)=F_j(V)\oplus F_j(W).
\]
Then,
$\theta$ takes $F_j(E)$ to $F_{j+1}(E)\otimes K_X$.

With respect to the filtration $\mathbf{F}(E)$,
the associated graded Higgs bundle is 
\begin{equation}\label{Graded1}
(E_0=M\oplus K_X^{n-1}\oplus K_X^{n-2}\oplus \cdots \oplus K_X^{2-n}\oplus K_X^{1-n}\oplus M^{-1}, \quad\theta_0=\begin{pmatrix}0&&&&&\\\mu&0&&&&\\&1&0&&&\\&&\ddots&\ddots&\\&&&1&0&\\&&&&\mu&0\end{pmatrix}).
\end{equation}

\begin{prop}\label{InitialCondition1}
Let $X$ be a non-compact hyperbolic Riemann surface.  Suppose $\mu\neq 0$. Suppose there exists a diagonal harmonic metric $h_1$ on $(E_0,\theta_0)$, compatible with $SO(n,n+1)$-structure. Then there exists a harmonic metric $h$ on $(E,\theta)$ which is compatible with $SO(n,n+1)$-structure and weakly dominates $h_1.$
\end{prop}
\pf
Let $X_i$ $(i=1,2,\cdots)$ be
a smooth exhaustion family of $X$.
Let $h^{(i)}$ be the harmonic metrics of $(E,\theta)|_{X_i}$
such that $h^{(i)}|_{\partial X_i}=h_0|_{\partial X_i}$.
Note that $h_0=h_0|_V\oplus h_0|_W$,
where $h_0|_V, h_0|_W$ are compatible with $Q_V, Q_W$ respectively. 
By Lemma \ref{Symmetry1},
$h^{(i)}$ is compatible with $SO(n,n+1)$-structure.
By Theorem \ref{Convergence}, $h^{(i)}$ has a convergence subsequence and has a smooth limit harmonic metric $h$. As a result, $h=h|_V\oplus h|_W,$ where $h|_V, h|_W$ are compatible with $Q_V, Q_W$ respectively.
\hfill\qed

\begin{thm}\label{ExistenceSO}
Suppose $X$ is the unit disk. Suppose there exists a flat Hermitian metric $h_M$ on $M$ and $\mu\in H^0(X, M^{-1}K_X^n)$ satisfies $h_M^{-1}g_X^{-n}(\mu,\mu)\in \mathcal A$ and not constantly $0$, then there exists a harmonic metric $h$ on $(E,\theta)$, compatible with $SO(n,n+1)$-structure.
\end{thm}
\pf
Consider
\begin{equation}
(E_1=M\oplus K_X^{n-1}\oplus K_X^{n-2}\oplus \cdots \oplus K_X^{2-n}\oplus K_X^{1-n}\oplus M^{-1}, \quad\theta_1=\begin{pmatrix}0&&&&&\\0&0&&&&\\&1&0&&&\\&&\ddots&\ddots&\\&&&1&0&\\&&&&0&0\end{pmatrix}).
\end{equation}
Let $$h_X=\oplus_{k=1}^{2n-1} a_{k,2n-1}g_X^{k-n}$$ be a diagonal metric on $$K_X^{n-1}\oplus K_X^{n-2}\oplus \cdots \oplus K_X^{2-n}\oplus K_X^{1-n},$$ where $a_{k,2n-1}$ is defined in Equation (\ref{a_kn}). 

Then $h_1=\diag(h_M, h_X,h_M^{-1})$ is a diagonal harmonic metric on $(E_1,\theta_1)$. We compare the Higgs bundle $(E_1,\theta_1)$ with $(E_0,\theta_0)$. 
It follows from Theorem \ref{NilpotentHiggsBundles} and $h_M^{-1}g_X^{-n}(\mu,\mu)\in \mathcal A$ that there exists a diagonal harmonic metric $h_0$ on the Higgs bundle $(E_0,\theta_0)$. Similar to the argument in Proposition \ref{InitialCondition1} and the proof of Theorem \ref{NilpotentHiggsBundles} and Proposition \ref{BasicExistence}, one can impose that $h_0$ is compatible with $SO(n,n+1)$-structure. Then the statement follows from Proposition \ref{InitialCondition1}.
\hfill\qed

\subsubsection{The generically regular semisimple case}
\label{subsection;23.6.16.3}


In this subsection,
we use the notation $(E,\theta(\vecq))$
to denote the Higgs bundle
associated with
$SO(n,n+1)$-Higgs bundle
$((V,Q_V),(W,Q_V),\eta_{\mu,\nu}(\vecq))$
in (\ref{OriginalBundle1})
to emphasize the dependence on $\vecq$.
According to Theorem \ref{GenericExistence},
if $((V,Q_V),(W,Q_V),\eta_{\mu,\nu}(\vecq))$
is generically regular semisimple,
then $(E,\theta(\vecq))$
has a harmonic metric compatible with
$SO(n,n+1)$-structure.
Let us mention some examples.

The following lemma is obvious.
\begin{lem}\label{GenericCondition}
If $\vecq=\mathbf{0}=(0,\ldots,0)$,
then
$\eta_{\mu,\nu}(\mathbf{0})^{\dagger}\circ \eta_{\mu,\nu}(\mathbf{0})
\in \End(V)\otimes K_X^2$ is induced by the identity morphisms
$K_X^{n+1-2i}\cong K_X^{n-1-2i}\otimes K_X^2$ $(i=1,\cdots, n-1)$
and $2\mu\nu:K_X^{-n+1}\rightarrow K_X^{n-1}\otimes K_X^2$.
Therefore,
if $\mu\nu$ is not constantly $0$, $((V,Q_V),(W,Q_W),\eta_{\mu,\nu}(\mathbf{0}))$ is generically regular semisimple.
\hfill\qed
\end{lem}

We obtain the following corollary from Lemma \ref{GenericCondition}
and Theorem \ref{GenericExistence}.
\begin{cor}
If $\vecq=0$ and if $\mu\nu$ is not constantly $0$, then there exists a harmonic metric $h$ of $(E,\theta(\mathbf{0}))$ compatible with $SO(n,n+1)$-structure. 
\hfill\qed
\end{cor}

Let us consider the case $X=\mathbb C$.
Let $M=\mathcal O_{\mathbb C}$.
Let $\mu_0$ and $\nu_0$ be non-zero polynomials.
We set $\mu=\mu_0dz^n$ and $\nu=\nu_0dz^n$.
For a positive integer $N$,
we set $\mathcal P_N=\{g(z)\in \mathbb C[z]\,|\,\deg g\leq N\}$.
We consider the following affine space.
$$\mathcal Q_N=
\{(g_1(z)dz^2,g_2(z)dz^4,\cdots, g_{n-1}(z)dz^{2n-2})
\,|\,g_{i}\in \mathcal P_N\}.
$$
\begin{prop}
There exists a non-empty Zariski open subset
$\mathcal U\subset \mathcal Q_N$
such that for any $\vecq\in \mathcal Q_N$
the associated $SO(n,n+1)$-Higgs bundle
$((V,Q_V),(W,Q_W),\eta_{\mu,\nu}(\vecq))$
is generically regular semisimple.
As a result,
for any $\vecq\in \mathcal U$,
the Higgs bundle
$(E,\theta(\vecq))$ on $\mathbb C$ has a harmonic metric compatible with $SO(n,n+1)$-structure. 
\end{prop}
\pf 
We obtain the first claim from Lemma \ref{GenericCondition}
which says $\mathbf{0}\in \mathcal U$ under the assumption $\mu\nu$
is not constantly $0$.
The second claim follows from Theorem \ref{GenericExistence}.
\hfill\qed

\section{$Sp(4,\mathbb R)$-Higgs bundles}\label{Sp(4,R)}

In this section, we discuss the existence of harmonic metrics on
$Sp(2n,\mathbb R)$-Higgs bundles over non-compact Riemann surfaces
by using the techniques developed in this paper and
our previous paper \cite{LiMochizukiGeneric}.
We are mainly interested in the case $n=2$.
\begin{df}
\begin{itemize}
 \item An $Sp(2n,\mathbb R)$-Higgs bundle over a Riemann surface $X$
       is determined by $(V, \gamma, \beta)$,
       where $V$ is a rank $n$ vector bundle,
       $\gamma\in H^0(X,S^2V^{\lor}\otimes K_X)$
       and $\beta\in H^0(X,S^2V\otimes K_X)$.
\item The associated $SL(2n,\mathbb C)$-Higgs bundle is $(E=V\oplus V^{\lor}, \theta=\begin{pmatrix}0&\beta\\
\gamma&0\end{pmatrix}).$
\item A harmonic metric $h$ on  $(E,\theta)$ is said to be compatible with $Sp(2n,\mathbb R)$-structure if $h=h|_V\oplus (h|_V)^{\lor}.$
\end{itemize}
\end{df}

The natural perfect pairing of $V$ and $V^{\lor}$
induces a non-degenerate symmetric pairing $Q_E$
of $E=V\oplus V^{\lor}$.
The Higgs field $\theta$ is self-adjoint with respect to $Q_E$.
If a harmonic $h$ of $(E,\theta)$ is compatible with
$\Sp(2n,\real)$-structure
then $h$ is compatible with $Q_E$.

\subsection{Dirichlet problem}

Let $Y\subset X$ be a relatively compact connected open subset
with smooth boundary $\del Y$.
Assume that $\del Y$ is non-empty.
Let $h_{\del Y}$ be any Hermitian metric of
$E_{|\del Y}$.
\begin{lem}\label{Symmetry2}
Let $h$ be a harmonic metric of $(E,\theta)$ such that $h|_{\del Y}=h_{\del Y}$. Suppose $h_{\del Y}$ is compatible with $Sp(2n,\mathbb R)$-structure. Then $h$ is compatible with $Sp(2n,\mathbb R)$-structure.
\end{lem}
\pf 
First we show that $h=h|_V\oplus h|_{V^{\lor}}.$ There exists the automorphism $\varphi=1_V\oplus (-1_{V^{\lor}})$ on $E=V\oplus V^{\lor}$. Because $\varphi^*\theta=-\theta$, $\varphi^*(h)$ is also a harmonic metric of $(E,\theta)$. By the uniqueness of the solution for Dirichlet problem for harmonic metric, $\varphi^*(h)|_{\partial Y}=h_{\partial Y}$, we obtain $\varphi^*(h)=h$. It means that $h$ is the direct sum of the Hermitian metrics of $V$ and $V^{\lor}$. 

Next we show that $h|_{V^{\lor}}=(h|_V)^{\lor}$. The metric $h$ induces the harmonic metric $h^{\lor}=(h|_V)^{\lor}\oplus (h|_{V^{\lor}})^{\lor}$ on $(E^{\lor}=V^{\lor}\oplus V, \theta^{\lor}=\begin{pmatrix}0&\gamma\\\beta&0\end{pmatrix}).$ Re-ordering $V$ and $V^{\lor}$, we have $(h|_{V^{\lor}})^{\lor}\oplus (h|_V)^{\lor}$ is a harmonic metric on $(E,\theta)$. Note that $\big((h|_{V^{\lor}})^{\lor}\oplus (h|_V)^{\lor}\big)|_{\del Y}=h_{\del Y}$. By the uniqueness of solutions of Dirichlet problem for harmonic metrics, we obtain $(h|_{V^{\lor}})^{\lor}\oplus (h|_V)^{\lor}=h|_V\oplus h|_{V^{\lor}}$. Thus, $h|_{V^{\lor}}=(h|_V)^{\lor}$.
\hfill\qed

\subsection{The generically regular semisimple case}

Let $(V,\gamma,\beta)$ be an $Sp(2n,\real)$-Higgs bundle on $X$.
Let $(E,\theta)$ denote the associated $SL(2n,\cnum)$-Higgs bundle.
We obtain
$\beta\circ\gamma\in\End(V)\otimes K_X^2$.
The spectral curve
$\Sigma_{\beta\circ\gamma}\subset (T^{\ast}X)^{\otimes 2}$
of $\beta\circ\gamma$ is defined as usual.
We obtain the finite map
$\pi:\Sigma_{\beta\circ\gamma}\to X$.

\begin{df}
$(V,\gamma,\beta)$ is called generically regular semisimple
if there exists $P\in X$
such that $|\pi^{-1}(P)|=n$
and $0\not\in\pi^{-1}(P)$. 
\hfill\qed 
\end{df}

The following theorem says $(E,\theta)$
has a harmonic metric compatible with $Sp(2n,\real)$-structure
in most cases.
See \S\ref{subsection;23.6.16.20}
for examples in the Gothen section.

\begin{thm}
\label{thm;23.6.16.30}
Suppose $X$ is a general non-compact Riemann surface.
If $(V,\gamma,\beta)$ is generically regular semisimple,
there exists a harmonic metric $h$ of
$(E,\theta)$ compatible with $Sp(2n,\mathbb R)$-structure.
\end{thm}
\pf
By Corollary \ref{Equivalence},
$(E,\theta)$ is generically regular semisimple.
It is standard to obtain the claim of Theorem \ref{thm;23.6.16.30}
by using Lemma \ref{Symmetry2},
\cite[Proposition 2.37]{LiMochizukiGeneric}
and Proposition \ref{prop;20.5.29.1}.
(See the proof of Theorem \ref{GenericExistence}.)
\hfill\qed

\begin{cor}
\label{cor;23.6.16.40}
Suppose $X$ is a general non-compact Riemann surface.
If
$\bigl(\tr(\beta\gamma)\bigr)^2-4\det\beta\cdot\det\gamma$
and $\det(\beta)\det(\gamma)$ are 
not constantly $0$,
there exists a harmonic metric $h$
of $(E,\theta)$ compatible with $Sp(4,\mathbb R)$-structure.
\end{cor}
\pf
If
$\bigl(\tr(\beta\gamma)\bigr)^2-4\det\beta\cdot\det\gamma$
and $\det(\beta)\det(\gamma)$
are not constantly $0$,
$(V,\beta,\gamma)$ is generically regular semisimple.
Hence, the claim follows from Theorem \ref{thm;23.6.16.30}.
\hfill\qed


\subsection{Gothen section}
Given a holomorphic line bundle $N$ on $X$ and
 $$\mu\in H^0(X, N^{-2}K_X^3),
 \quad
 \nu\in H^0(X, N^2K_X),
 \quad
 q_2\in H^0(X, K_X^2),$$
 one can construct 
a $Sp(4,\mathbb R)$-Higgs bundle
$(V,\gamma,\beta)$
 as follows:
\[
V=N\oplus N^{-1}K_X,
\quad
\gamma=\begin{pmatrix}0&1\\1&0\end{pmatrix},
\quad
\beta=\begin{pmatrix}\nu&q_2\\q_2&\mu\end{pmatrix}.
\] 
The associated $SL(4,\mathbb C)$-Higgs bundle is
\begin{equation}\label{OriginalBundle2}
(E=N\oplus N^{-1}K_X\oplus N^{-1}\oplus NK_X^{-1},\quad\theta=\begin{pmatrix}0&0&\nu&q_2\\
0&0&q_2&\mu\\
0&1&0&0\\
1&0&0&0
\end{pmatrix}).
\end{equation}

When $X$ is a compact Riemann surface of genus at least two, for each integer $d\in (g-1, 3g-3],$ there is a component $X_d$ (see \cite{Gothen}, \cite[Proposition 3.23]{BradlowOscarGothen}), called Gothen component, of the moduli space of  $Sp(4,\mathbb R)$-Higgs bundles formed by the above Higgs bundles determined by $(N,\mu,\nu,q_2)$ where $\deg(N)=d$ and $\mu\neq 0$. In particular, when $d=3(g-1), $ it coincides with the Hitchin component for $Sp(4,\mathbb R).$ Such Higgs bundles are maximal and correspond to maximal $Sp(4,\mathbb R)$ representations. 

We call the above Higgs bundles are in the Gothen section. We are going to discuss the existence of harmonic metrics of Higgs bundles in the Gothen section over non-compact Riemann surfaces.


\subsubsection{The generically regular semisimple case}
\label{subsection;23.6.16.20}

\begin{prop}
Suppose $X$ is a non-compact Riemann surface.
If $\mu\nu$ and $\mu\nu-q_2^2$
are not constantly $0$,
there exists a harmonic metric $h$ of $(E,\theta)$
compatible with $Sp(4,\mathbb R)$-structure.
\end{prop}
\pf
Because
$\det(\beta\gamma)=q_2^2-\mu\nu$
and
$(\tr\beta\gamma)^2-4\det(\beta\gamma)
=(2q_2)^2-4(q_2^2-\mu\nu)=4\mu\nu$,
we obtain the claim 
from Corollary \ref{cor;23.6.16.40}.
\hfill\qed


\subsubsection{The case $(\mu,\nu)=(0,0)$}
\begin{prop}
Suppose $X$ is a non-compact Riemann surface. Suppose in addition $q_2\neq 0$ when $X$ is parabolic. If $(\mu,\nu)=(0,0)$, then there exists a harmonic metric $h$ of $(E,\theta)$ compatible with $Sp(4,\mathbb R)$-structure.
\end{prop}
\pf
For $(\mu,\nu)=(0,0)$, the Higgs bundle $$(E,\theta)=\big(N\oplus NK_X^{-1},\begin{pmatrix}0&q_2\\1&0\end{pmatrix}\big)\oplus \big(N^{-1}K_X\oplus N^{-1},\begin{pmatrix}0&q_2\\1&0\end{pmatrix}\big).$$

Fix a square root line bundle $K_X^{\frac{1}{2}}$ of $K_X$.
Let $L=NK_X^{-\frac{1}{2}}$ and let $h_L$ be a flat Hermitian metric on $L$.
Let $\diag(h_0, h_0^{-1})$ be a harmonic metric on $\big (K_X^{\frac{1}{2}}\oplus K_X^{-\frac{1}{2}}, \begin{pmatrix}0&q_2\\1&0\end{pmatrix}\big).$
Then $$\diag(h_L\otimes h_0, h_L\otimes h_0^{-1})\oplus \diag(h_L^{-1}\otimes h_0, h_L^{-1}\otimes h_0^{-1})$$ is a harmonic metric of $(E,\theta)$ compatible with $Sp(4,\mathbb R)$-structure.
\hfill\qed

\subsubsection{The case $\mu\neq 0$}
Set $$F_1=N, \quad F_2=N\oplus NK_X^{-1}, \quad F_3=N\oplus NK_X^{-1}\oplus N^{-1}K_X,\quad  F_4=E.$$ Then $\mathbf{F}=\{F_1\subset F_2\subset F_3\subset F_4\}$ is a full holomorphic filtration of $E$ and $\theta$ takes $F_i$ to $F_{i+1}\otimes K.$ And the graded Higgs bundle is 
\begin{equation}\label{Graded2}
(E_0=N\oplus NK_X^{-1}\oplus N^{-1}K_X\oplus N^{-1},\quad \theta_0=\begin{pmatrix}0&0&0&0\\
1&0&0&0\\
0&\mu&0&0\\
0&0&1&0
\end{pmatrix}).
\end{equation}
\begin{prop}\label{InitialCondition2}
Let $X$ be a non-compact hyperbolic Riemann surface. Suppose $\mu\neq 0$. Suppose there exists a diagonal harmonic metric $h_1$ on $(E_0,\theta_0)$, compatible with $Sp(4,\mathbb R)$-structure. Then there exists a harmonic metric $h$ on $(E,\theta)$ which is compatible with $Sp(4,\mathbb R)$-structure and weakly dominates $h_1$ with respect to $\mathbf F$.
\end{prop}
\pf
Let $X_i$ $(i=1,2,\cdots)$ be
a smooth exhaustion family of $X$.
Let $h^{(i)}$ be the harmonic metrics of $(E,\theta)|_{X_i}$
such that $h^{(i)}|_{\partial X_i}=h_1|_{\partial X_i}$.
Note that $h_1=h_1|_V\oplus (h_1|_V)^{\lor}$ for $V=N\oplus N^{-1}K_X$. 
By Lemma \ref{Symmetry2}, we also have $h^{(i)}=h^{(i)}|_V\oplus (h^{(i)}|_V)^{\lor}$.
By Theorem \ref{Convergence}, $h^{(i)}$ has a convergence subsequence and limits to a harmonic metric $h$ which weakly dominates $h_1$ with respect to $\mathbf F$. As a result, $h=h|_V\oplus (h|_V)^{\lor}$.
\hfill\qed
\vspace{0.1cm}

Fix a square root $K_X^{\frac{1}{2}}$ of $K_X$. 
\begin{thm}
Suppose $X$ be the unit disk. Let $L=NK_X^{-\frac{1}{2}}$. Suppose there exists a flat Hermitian metric $h_L$ on $L$. If $\mu\in H^0(X, L^{-2}K_X^2)$ satisfies $h_L^{-2}g_X^{-2}(\mu,\mu)\in \mathcal A$ and not constantly $0$, then there exists a harmonic metric $h$ on $(E,\theta)$, compatible with $Sp(4,\mathbb R)$-structure.
\end{thm}
\pf
The Higgs bundle $(E_0,\theta_0)$ becomes 
\begin{equation}\label{Graded3}
(E_0=LK_X^{\frac{1}{2}}\oplus LK_X^{-\frac{1}{2}}\oplus L^{-1}K_X^{\frac{1}{2}}\oplus L^{-1}K_X^{-\frac{1}{2}},\quad \theta_0=\begin{pmatrix}0&0&0&0\\
1&0&0&0\\
0&\mu&0&0\\
0&0&1&0
\end{pmatrix}),
\end{equation} where $\mu\in H^0(X, L^{-2}K_X^2)$.

Consider Higgs bundle 
\begin{equation}\label{Graded4}
(E_1=LK_X^{\frac{1}{2}}\oplus LK_X^{-\frac{1}{2}}\oplus L^{-1}K_X^{\frac{1}{2}}\oplus L^{-1}K_X^{-\frac{1}{2}},\quad \theta_1=\begin{pmatrix}0&0&0&0\\
1&0&0&0\\
0&0&0&0\\
0&0&1&0
\end{pmatrix}),
\end{equation}
equipped with a diagonal harmonic metric $$h_1=\diag(\sqrt{2}h_Lg_X^{-\frac{1}{2}}, \frac{1}{\sqrt{2}}h_Lg_X^{\frac{1}{2}}, \sqrt{2}h_L^{-1}g_X^{-\frac{1}{2}},\frac{1}{\sqrt{2}}h_L^{-1}g_X^{\frac{1}{2}}).$$ We compare the Higgs bundle $(E_1,\theta_1)$ with $(E_0,\theta_0)$. 
It follows from Theorem \ref{NilpotentHiggsBundles} and $h_L^{-2}g_X^{-2}(\mu,\mu)\in \mathcal A$ that there exists a diagonal harmonic metric $h_0$ on the Higgs bundle $(E_0,\theta_0)$. 

By Lemma \ref{Symmetry2} and the proof of Theorem \ref{NilpotentHiggsBundles} and Proposition \ref{BasicExistence}, one can impose that $h_0$ is compatible with $Sp(4,\mathbb R)$-structure. 
\hfill\qed

\appendix
\section{Discussions on Green functions}

In this section, let $X$ denote the unit disk.
Recall that the Green function on $X$ is 
$$G(z,\xi)=\log\Big|\frac{1-z\bar \xi}{z-\xi}\Big|.$$

\begin{lem}\label{PropertiesGreenFunction}
1. For a nonnegative function $f$ on $X,$ the followings are equivalent:
\begin{itemize}
    \item 
$f\in \mathcal A.$
\item for some $z$,  $\int_Xf(\xi)G(z,\xi)d\sigma_\xi< \infty$. 
\item for all $z$,  $\int_Xf(\xi)G(z,\xi)d\sigma_\xi< \infty$.
\end{itemize}
2. For a bounded function $u$ and a nonnegative function $f$ on $X$, $f\in \mathcal A^b$ if and only $e^u\cdot f\in \mathcal A^b.$\\
3. For a nonnegative function $f$ on $X,$ if $f(z)\leq M(1-|z|^2)^p$ for a constant $p>-2,$ then $f\in \mathcal A^b.$\\ 
4. For a constant $p\leq -2,$ then $(1-|z|^2)^p\notin \mathcal A.$
\end{lem}
\pf
For Part 1, for $0\leq r_0<s<1$, 
\[G(r_0,s)=\log\frac{|1-r_0s|}{|r_0-s|}=\log(1+\frac{(1-s)(1+r_0)}{s-r_0})\]
\[G(-r_0,s)=\log\frac{|1+r_0s|}{|r_0+s|}=\log(1+\frac{(1-s)(1-r_0)}{s+r_0}).\]

So there exist positive constants $C_1, C_2,C_3, C_4$ such that 
\[C_1\frac{(1-s)(1+r_0)}{s-r_0}\leq G(r_0,s)\leq C_2\frac{(1-s)(1+r_0)}{s-r_0},\]
\[C_3\frac{(1-s)(1+r_0)}{s-r_0}\leq G(-r_0,s)\leq C_4\frac{(1-s)(1+r_0)}{s-r_0}.\]
Together with the inequality $G(-|z|,|\xi|)\leq G(z,\xi)\leq G(|z|,|\xi|)$ (see \cite[Lemma 4.5.7]{Ransford}), the statements in Part 1 follows.

Part 2 is clear.

For Part 3 and Part 4, for $z\neq 0,$ 
\begin{eqnarray*}
&& \int _XG(z,\xi)(1-|z|^2)^pd\sigma_\xi \\
&&\text{write $\xi=r^{i\theta}$}\\
&=& \frac{1}{2\pi}\int_0^1(\int_0^{2\pi}\log\frac{|1-z\cdot re^{-i\theta}|}{ |z-re^{i\theta}|}d\theta)(1-r^2)^prdr \\
&=& \frac{1}{2\pi}\int_0^1(\int_0^{2\pi}(\log|z|+\log|\frac{1}{\bar z}-re^{i\theta}|-\log |z-re^{i\theta}|)d\theta)(1-r^2)^prdr \\
&&\text{Using Lemma \ref{BasicCalculation}}\\
&=& \int_0^1(\log|z|+\log(\max\{\frac{1}{|z|}, r\})-\log(\max\{|z|, r\}))(1-r^2)^prdr \\
&=& \int_0^1(\log(\max\{1, r|z|\})-\log(\max\{|z|, r\}))(1-r^2)^prdr \\
&=&-\int_0^1\log(\max\{|z|, r\})(1-r^2)^prdr 
\end{eqnarray*}
So \begin{eqnarray*}
&& \int _XG(z,\xi)(1-|z|^2)^pd\sigma_\xi\\
&\leq&-\int_0^1\log(\max\{|z|, r\})(1-r)^prdr \\
&\leq&-\int_0^1\log r(1-r)^prdr 
\end{eqnarray*}
which is integrable for $p>-2$. Together with Part 1, the integral is finite when $z=0$. 
Therefore $$\sup_{z\in X}\int _XG(z,\xi)(1-|\xi|^2)^pd\sigma< \infty.$$ So Part 3 follows. 

 \begin{eqnarray*}
&& \int _XG(z,\xi)(1-|z|^2)^pd\sigma_\xi\\
&\geq&-\int_0^1\log(\max\{|z|, r\})(1-r)^prdr \\
&\geq&-[\int_{|z|}^1\log r(1-r)^prdr ],
\end{eqnarray*} which is infinity for $p\leq -2.$
So Part 4 follows.
\hfill\qed

\begin{lem}\label{BasicCalculation}
For all $z\in \mathbb C$ any $r>0$, 
\[\frac{1}{2\pi}\int_0^{2\pi}\log|z-re^{i\theta}|d\theta=\log(\max\{|z|,r\}).\]
\end{lem}
\pf
The function $\log |z|$ is harmonic on $\mathbb C^*.$ If $r<|z|,$ then $\overline{\mathbb D(z,r)}\subset \mathbb C^*.$ By the mean value property of harmonic functions, we have
\[\frac{1}{2\pi}\int_0^{2\pi}\log|z-re^{i\theta}|d\theta=\log|z|.\]

Now suppose $r>|z|.$ Write $\log|z-re^{i\theta}|=\log |z|+\log r+\log |\frac{1}{z}-\frac{1}{r}e^{i\theta}|.$ Then we have 
\[\frac{1}{2\pi}\int_0^{2\pi}\log|z-re^{i\theta}|d\theta=\log|z|+\log r+\frac{1}{2\pi}\int_0^{2\pi}\log |\frac{1}{z}-\frac{1}{r}e^{i\theta}|.\]
Note that $|\frac{1}{z}|<\frac{1}{r},$ then $\frac{1}{2\pi}\int_0^{2\pi}\log |\frac{1}{z}-\frac{1}{r}e^{i\theta}|=\log\frac{1}{z}$. Therefore \[\frac{1}{2\pi}\int_0^{2\pi}\log|z-re^{i\theta}|d\theta=\log|z|+\log r+\log\frac{1}{|z|}=\log r.\]

By the continuity, we also obtain that when $|z|=r$, the integral is zero.

\section{Various expressions of Higgs bundles in the Hitchin section}\label{CharacterizationHitchinSection}
Given a $SL(n,\mathbb C)$-Higgs bundle $(E,\theta)$ satisfies that $E$ admits a full holomorphic filtration $\mathbf{F}=\{F_1\subset F_2\subset\cdots\subset F_n\}$ and $\theta(\vecq)$ takes $F_i$ to $F_{i+1}\otimes K_X$ and induces an isomorphism between $F_i/F_{i-1}\rightarrow F_{i+1}/F_i\otimes K$ for $i=1,\cdots, n-1$.  On a compact hyperbolic Riemann surface, there always exists a gauge transformation taking $(E,\theta)$ to $(\hyperk_{X,n}, \theta(\vecq))$ for some $\vecq$. One may refer \cite{CollierWentworth} for the proof. Therefore, Higgs bundles in the Hitchin section are characterized by the above filtration description.

Here we restrict to a more special description of $(E,\theta)$ and consider an arbitrary (possibly non-compact) Riemann surface $X$. Fix a choice of square root $K_X^{\frac{1}{2}}$ and let $L_i=(K_X^{\frac{1}{2}})^i$, for $1\leq i\leq n$. Let $E:=\hyperk_{X,n}=\oplus_{i=1}^n L_i$ and $\theta\in H^0(X, \End(E)\otimes K_X)$ be as follows: 
 Decompose $\End(E)=\oplus_{1\leq i,j\leq n}\Hom(L_i,L_j)$ and write $\theta=\theta_{ij}$ such that $\theta_{ij}\in \Hom(L_i,L_j)\otimes K_X$. Then $\theta_{jk}=0$ for all $j<k$ and $\theta_{j,j+1}\in H^0(X, \Hom(L_j,L_{j+1})\otimes K_X)\cong H^0(X, \mathcal O)$ is a nonzero constant for each $j=1,2, \cdots, n-1$. Moreover, $\Tr(\theta)=0$.
 
For $j\geq 0$, under the isomorphism $\Hom(L_k,L_{k-j+1})\otimes K_X\cong K_X^j$, let $q_j=\frac{1}{n+1-j}\sum_{k=j}^n\theta_{k,k-j+1}\in H^0(X,K_X^j)$. It follows from $\tr(\theta)=0$ that $q_1=0$. For such Higgs bundle $(\hyperk_{X,n},\theta)$, we can construct an explicit holomorphic gauge transformation taking it to $(\hyperk_{X,n}, \theta(\vecq)). $ The construction is similar to the one in \cite[Proposition 3.11]{CollierWentworth} except here we do not need to use the invertibility of an elliptic operator. We include the construction in the following for its own interest. 
\begin{prop}
There exists a holomorphic gauge transformation $g$ taking $(\hyperk_{X,n},\theta)$ to $ (\hyperk_{X,n}, \theta(\vecq)). $ That is, $$g\circ\theta \circ g^{-1}=\theta(\vecq)=\begin{pmatrix}0&q_2&q_3&q_4&\cdots&q_n\\
1&0&q_2&q_3&\ddots&\vdots\\
&1&0&q_2&\ddots&\vdots\\
&&\ddots&\ddots&\ddots&q_3\\
&&&\ddots&\ddots&q_2\\
&&&&1&0
\end{pmatrix}.$$
\end{prop}
\pf Let $r_k=\theta_{k,k+1}$. 
One could choose the gauge transformation $g$ which is formed by $\prod_{l=1}^{k-1}r_l^{-1}$ on $\Hom(L_k,L_k)\cong \mathcal O$, for $k=1,2,\cdots,n$. Then $g\theta g^{-1}$ satisfies that $(g\theta g^{-1})_{j,j+1}=1.$ So we can always assume $\theta_{j,j+1}=1$ for $1\leq j\leq n-1$. 

Set $\End(E)_j=\oplus_{i-k=j}\Hom(L_i,L_k)$. For $f\in \Omega^i(X,\End(E)_k), g\in \Omega^j(X,\End(E)_l)$, the composition $f\wedge g\in \Omega^{k+l}(X,\End(E)_{k+l})$. For $f\in \Omega^k(X,\End(E))$, one may decompose $f=\oplus_{j=1-n}^{n-1}f_j$ for $f_j\in \Omega^k(X,\End(E)_j)$. Denote by $f_{\geq j}$ the sum of $f_k$'s for $k\geq j$. So $\theta=\theta_{\geq -1}$. 

We are going to prove the proposition by induction on $j$, starting from $0$. 
For $j=0$, $\theta_{k,k-j+1}=1$ for $1\leq k\leq n-1$.

Suppose for $i\leq j_0$, $\theta_{k,k-i+1}=q_i$ for $i\leq k\leq n$. 
We will find a holomorphic gauge transformation $g$ such that $(g\circ\theta\circ g^{-1})_{k,k-i+1}=q_{i}$ for each $i\leq k\leq n$ and $i\leq j_0+1$. 
Choose $g=I+f$, where $I$ is the identity matrix and $f=f_{j_0+1}$. 
So $g^{-1}=I-f+h$ for some $h$ satisfying $h=h_{\geq j_0+2}$. 
Then $$g\circ \theta\circ g^{-1}=(I+f)\theta(I-f+h).$$ Note that $\theta=\theta_{\geq -1}, f=f_{j_0+1}, h=h_{\geq j_0+2}$. 
So for $i\leq j_0$, 
$$(g\theta g^{-1})_{i-1}=\theta_{i-1}. $$ By assumption, $(g\circ\theta \circ g^{-1})_{k,k-i+1}=q_i$ for $i\leq k\leq n$.

 For $i=j_0+1$,
$$(g\circ \theta\circ g^{-1})_{j_0}=\theta_{j_0}+f_{j_0+1}\circ\theta_{-1}-\theta_{-1}\circ f_{j_0+1}.$$
Then $(g\circ\theta\circ g^{-1})_{k,k-j_0}=\theta_{k,k-j_0}+f_{k+1,k-j_0}-f_{k,k-j_0-1}$. 
Set $$f_{k,k-j_0-1}=-\sum_{l=j_0+1}^{k-1}\theta_{l,l-j_0}+(k-j_0-1)q_{j_0+1}.$$ 
We obtain that $(g\circ \theta \circ g^{-1})_{k,k-j_0}=q_{j_0+1}$ for $j_0+1\leq k\leq n$. 
Clearly $g$ is holomorphic, following from $\theta_{ij}$ are holomorphic. 
\hfill\qed

\bibliographystyle{amsalpha}
\bibliography{bib}

\end{document}